\documentclass[preprint,12pt]{elsarticle}

\usepackage{lipsum}
\makeatletter
\def\ps@pprintTitle{%
	\let\@oddhead\@empty
	\let\@evenhead\@empty
	\def\@oddfoot{}%
	\let\@evenfoot\@oddfoot}
\makeatother

\usepackage[margin=1in]{geometry}

\usepackage{graphicx}
\usepackage{amssymb}

\usepackage{lineno}
\usepackage{amsmath}
\usepackage{amssymb}
\usepackage{amsfonts}

\usepackage[section]{placeins}

\usepackage{float}

\usepackage{graphics}
\usepackage{epsfig}

\usepackage{color}
\usepackage{hyperref}
\usepackage{comment}
\usepackage{multirow}
\usepackage{array}
\usepackage[dvipsnames]{xcolor}

\usepackage{multicol,lipsum}
\usepackage{subcaption}
\usepackage{soul}
\usepackage[utf8]{inputenc}
\usepackage[font=normal]{caption}
\usepackage[version=4]{mhchem}
\usepackage{siunitx}
\usepackage{longtable,tabularx}
\usepackage{makecell}
\usepackage{cancel}

\usepackage{ulem}

\makeatletter
\DeclareRobustCommand{\qed}{%
	\ifmmode 
	\else \leavevmode\unskip\penalty9999
	\hbox{}\nobreak\hfill
	\fi
	\quad\hbox{\qedsymbol}}
\newcommand{\openbox}{\leavevmode
	\hbox to.77778em{%
		\hfil\vrule
		\vbox to.675em{\hrule width.6em\vfil\hrule}%
		\vrule\hfil}}
\newcommand{\qedsymbol}{\openbox}

\newcommand{\proofname}{Proof}
\makeatother

\usepackage{algorithm}
\usepackage{algpseudocode}

\algtext*{EndFor}
\algtext*{EndIf}
\algtext*{EndFunction}
\algtext*{EndWhile}

\algrenewcommand\algorithmicforall{\textbf{foreach}}
\algrenewcommand\algorithmicindent{.8em}

\journal{}

\begin{document}
	
\begin{frontmatter}

\title{Optimization under Uncertainty of a Chain of Nonlinear Resonators using a Field Representation}

\author{Seyed Saeed Ahmadisoleymani}
\ead{ssaeeda@email.arizona.edu}

\author{Samy Missoum\corref{cor}}
\ead{smissoum@arizona.edu}

\address{Department of Aerospace and Mechanical Engineering, University of Arizona, Tucson, Arizona}

\cortext[cor]{Corresponding author}

\begin{abstract}

Chains of resonators in the form of spring-mass systems have long been known to exhibiting interesting properties such as band gaps.  Such features can be leveraged to manipulate the propagation of waves such as the filtering of specific frequencies and, more generally, mitigate vibrations and impact. Adding nonlinearities to the system can also provide further avenues to manipulate the propagation of waves in the chain and enhance its performance.
This work proposes to optimally design such a chain of resonators to mitigate vibrations in a robust manner by accounting for various sources of design uncertainties (e.g., nonlinear stiffness) and aleatory uncertainties (e.g., loading). The stochastic optimization algorithm is tailored to account for discontinuities in the chain response due to the presence of nonlinearities. In addition, a field formulation is used to define the properties of the resonators along the chain and reduce the dimensionality of the optimization problem. It is shown that the combination of the stochastic optimization algorithm and the field representation leads to robust designs that could not be achieved with optimal properties constant over the chain.

\end{abstract}

\end{frontmatter}

\section{Introduction}
\label{sec:introduction}
The design of metamaterials has been the object of many studies because interesting dynamic behavior, sometimes non-existent in nature, can be achieved. For instance, a basic phenomenon extensively studied is the appearance of band gaps, which are frequency ranges within which waves cannot propagate \cite{sigalas1992elastic, kushwaha1994theory}. There is a large body of literature related to various metamaterial designs exploiting band gaps \cite{pai2010metamaterial,casadei2012piezoelectric,beli2018wave}. Techniques such as topology optimization applied to the design of unit cells in periodic metamaterials have also been used extensively to tailor such band gaps \cite{sigmund2003systematic, halkjaer2006maximizing,dahl2008topology, vatanabe2014maximizing, liu2016systematic}.

Beyond the well-known band gaps, it has long been understood that nonlinearities at the unit cell level can substantially increase the range of possible wave manipulations. There again, there have been many studies exploiting nonlinearities in various forms (e.g., stiffness, constitutive law) to achieve specific wave characteristics  \cite{narisetti2010perturbation, manktelow2013topology, manimala2016numerical,daraio2006tunability}. Leveraging band gaps and other effects (e.g., targeted energy transfer) due to local nonlinearities can markedly improve the capabilities of metamaterials. Among all the possible applications, there is a promising avenue in the area of vibration and shock mitigation, which is the focus application of this study.

The design optimization of metamaterials is tedious for several reasons. First, nonlinearities can lead to discontinuous responses. {These discontinuities demonstrate that the system's dynamic behavior is potentially highly sensitive to uncertainties, which must therefore be accounted for in the design process.} Examples of sources of uncertainties are design uncertainties (e.g., material properties), loading uncertainties, or initial conditions. Finally, the number of cells in a metamaterial is very large. Therefore, if a non-periodic structure is considered, the individual optimization of each unit is computationally very expensive. It is noteworthy that the non-smoothness of the problem is a major hurdle for topology optimization techniques, which are based on an efficient evaluation of the gradients using the adjoint method. This limitation is even more stringent because of the need to include uncertainties in the optimization problem. 

This work proposes a methodology to perform the design optimization under the uncertainty of a chain of nonlinear resonators, representative of a metamaterial. One of the key features of the stochastic optimization method lies in its ability to account for discontinuities, which are automatically detected through clustering \cite{boroson2017stochastic}. The clusters correspond to various levels of performance of the system and are associated with different regions of the search space. These regions are identified using a support vector machine (SVM) classifier \cite{christianini2000support},  which is locally refined using a dedicated adaptive sampling scheme \cite{lacaze2014generalized}. The performance in each region is approximated using  Kriging metamodels. The use of Kriging approximations makes it possible to efficiently propagate uncertainties to compute expected values based on Monte-Carlo simulations. The approach was originally developed for nonlinear energy sinks (NESs) which exhibit a discontinuous response in the form of an activation threshold \cite{boroson2017optimization, boroson2017stochastic}. The nonlinear resonators embedded in the metamaterial studied in this article might also exhibit discontinuous responses, thus requiring the need for a specific optimization approach. The stochastic optimization method is able to account for design uncertainties (e.g., stiffness properties) and ``aleatory" uncertainties (e.g., loading conditions). {Note that in this study, the wording ``aleatory" is used to describe random parameters representing loading conditions, and differentiate them from the random design variables. The proposed method} also addresses the difficulty due to the dimensionality of the optimization problem (i.e., a large number of unit cells) through the use of field representation. For instance, in order to reduce the dimensionality and make the problem scalable, properties such as the nonlinear stiffness are approximated using a spatial field over the whole chain. The field is governed by a handful of coefficients, thus substantially reducing the dimensionality of the problem. 

Among the possible problem formulations, the proposed approach is used to optimize the mean of performance metrics such as the RMS response in the case of vibration mitigation. In this study, the optimization of a one-dimensional chain of resonators is considered for demonstrative purposes. The chain  is made of units constituted of a main and an internal nonlinear resonators. This so-called ``mass-in-a-mass" design is originally representative of inclusions within a material. The nonlinearities considered in this article are in the form of a cubic stiffness property in the internal resonators.

This article is structured as follows. Section \ref{sec:chain_design} describes the chain of resonators and its properties, and provides the governing equations. The derivation of the dispersion analysis of the chain is provided to identify band gaps and understand how they vary as a function of the nonlinearity of the internal resonator. Section \ref{sec:optim_formul} describes the optimization problem formulation, the field implementation for dimensionality reduction, and the main elements of the stochastic optimization algorithm.  Results are provided in Section \ref{sec:results}, including three examples of stochastic optimization with harmonic loading.

\section{Chain of nonlinear resonators}	
\label{sec:chain_design}
\label{sec:chain_des}
\begin{figure}[!h]
	\centering
	\includegraphics[width=1\linewidth]{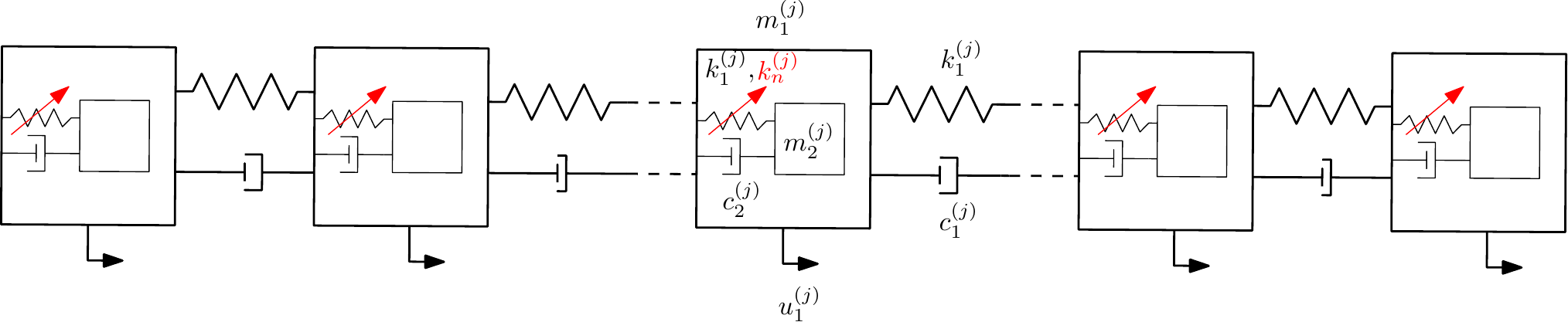}
	\caption{\small{Schematic representation of the "mass-in-a-mass" chain of resonators. The internal resonators are characterized by a linear and a cubic nonlinear stiffness.}}
	\label{fig:chain}
\end{figure}

In this study, a discrete one-dimensional chain of resonators is considered (Figure \ref{fig:chain}). Each unit of the chain contains an internal resonator thus representing a ``mass-in-a-mass" design \cite{huang2009negative}. The main resonators are connected through linear springs and linear dashpots. The internal resonators are connected to the main resonators through springs with linear and cubic stiffness properties and a linear dashpot. \\

\subsection{Governing equations}
The governing equations for the $j^{th}$ unit composed of a main and an internal resonators without external load are:

\begin{align}
\label{eqn:eqn_main}
& m_1^{(j)} \ddot{u}_1^{(j)} + k_1^{(j)}(u_1^{(j)}-u_1^{(j+1)}) + k_1^{(j-1)}(u_1^{(j)}-u_1^{(j-1)}) + k_2^{(j)}(u_1^{(j)}-u_2^{(j)}) + \dots \notag \\ & c_1^{(j)}(\dot{u}_1^{(j)}-\dot{u}_1^{(j+1)}) + c_1^{(j-1)}(\dot{u}_1^{(j)}-\dot{u}_1^{(j-1)})  + c_2^{(j)}(\dot{u}_1^{(j)}-\dot{u}_2^{(j)}) + k_n^{(j)}(u_1^{(j)}-u_2^{(j)})^3= 0 \notag \\
& m_2^{(j)} \ddot{u}_2^{(j)} + k_2^{(j)}(u_2^{(j)}-u_1^{(j)}) + c_2^{(j)}(\dot{u}_2^{(j)}-\dot{u}_1^{(j)}) + k_n^{(j)}(u_2^{(j)}-u_1^{(j)})^3 = 0 
\end{align}

\noindent where $u_1$ and $u_2$ are the displacements of the main and internal resonators, respectively. $k_1$ and $k_2$ are the linear stiffnesses. $c_1$ and $c_2$ are the viscous damping coefficients. $k_n$ is the nonlinear stiffness.

After non-dimensionalization, and removal of the unit index $j$ for readability, the equations become:

\begin{align}
\label{eqn:eqn_main}
& \bar{\omega}^2v''_1 + (2v_1-v_1^{(+1)}-v_1^{(-1)}) + \bar{\omega}\xi_1(2v'_1-{v'_1}^{(+1)} - {v_1'}^{(-1)}) + \alpha(v_1-v_2) + \dots \notag \\ & \bar{\omega}\xi_2(v'_1-v'_2) + \eta(v_1-v_2)^3 = 0 \notag \\
& \bar{\omega}^2\epsilon v''_2 + \alpha(v_2-v_1)+ \bar{\omega}\xi_2({v'_2}-{v'}_1) + \eta(v_2-v_1)^3  = 0 
\end{align}

where $'$ is the differentiation with respect to the non-dimensional time $\tau=\omega t$,  $\displaystyle{v=\frac{u}{L}}$ (where $L$ is the spatial distance between adjacent units), mass ratio $\epsilon=\displaystyle{\frac{m_2}{m_1}}$, stiffness ratio $\alpha=\displaystyle{\frac{k_2}{k_1}}$, frequency parameter $\bar{\omega}=\displaystyle{\frac{\omega}{\omega_1}}$ (where $\omega_1=\displaystyle{\sqrt{\frac{k_1}{m_1}}}$), damping terms $\xi_1=\displaystyle{\frac{c_1}{m_1\omega_1}}$ and $\xi_2=\displaystyle{\frac{c_2}{m_1\omega_1}}$, and nonlinear stiffness term $\eta=\displaystyle{\frac{k_nL^2}{m_1\omega_1^2}}$. Note that $L$ does not appear explicitly in the dimensional equations of motion.

\subsection{Band gaps}
\label{sec:HM}
In order to determine the frequency band gaps and assess the effect of nonlinearities, the dispersion behavior of the chain is studied using the harmonic balance method, as described in \cite{guskov2007multi, narisetti2012study}. This study is carried out without the presence of damping. The dimensional equations of motion for the $j^{th}$ unit can be written in matrix form as follows:
\begin{align}
\label{eqn:eqn_main2}
&\boldsymbol{M}\frac{d^2 \boldsymbol{u}^{(j)}}{dt^2}+\boldsymbol{K}^{(j-1)}\boldsymbol{u}^{(j-1)}+\boldsymbol{K}^{(j)}\boldsymbol{u}^{(j)}+\boldsymbol{K}^{(j+1)}\boldsymbol{u}^{(j+1)}+\boldsymbol{f}^{(j)}_{nl}=\boldsymbol{f}_{ext}
\end{align}

\noindent where	$\boldsymbol{M}=\begin{bmatrix} m_1 & 0\\ 0 & m_2 \end{bmatrix}$,
$\boldsymbol{u}^{(j)}=\begin{bmatrix} u_1^{(j)}\\u_2^{(j)}\end{bmatrix}$,
$\boldsymbol{K}^{(j-1)}=\begin{bmatrix} -k_1 & 0\\ 0 & 0 \end{bmatrix}$,\\
$\boldsymbol{K}^{(j+1)}=\begin{bmatrix} -k_1 & 0\\ 0 & 0 \end{bmatrix}$,
$\boldsymbol{K}^{(j)}=\begin{bmatrix} 2k_1+k_2 & -k_2\\ -k_2 & k_2 \end{bmatrix}$,
$\boldsymbol{f}^{(j)}_{nl}=\begin{bmatrix} k_n(u_1^{(j)}-u_1^{(j)})^3\\ -k_n(u_1^{(j)}-u_1^{(j)})^3\end{bmatrix}$  is the force term due to the nonlinear springs
and $\boldsymbol{f}_{ext}$ is the external load. \\

Based on the non-dimensional time $\tau=\omega t$, 	and setting the external load to zero, Eq. \ref{eqn:eqn_main2} can be condensed as:
\begin{align}
\label{eqn:eqn_main3}
&\omega^2\boldsymbol{M}\frac{d^2 \boldsymbol{u}^{(j)}}{d\tau^2}+ \boldsymbol{F_{int}}(\boldsymbol{u}^{(j)},\boldsymbol{u}^{(j-1)},\boldsymbol{u}^{(j+1)}) = 0
\end{align}
where,
\begin{align}
\label{eqn:finternal}
&\fontsize{11.4}{12} \boldsymbol{F_{int}}(\boldsymbol{u}^{(j)},\boldsymbol{u}^{(j-1)},\boldsymbol{u}^{(j+1)}) =\boldsymbol{K}^{(j-1)}\boldsymbol{u}^{(j-1)} + \boldsymbol{K}^{(j)}\boldsymbol{u}^{(j)} + \boldsymbol{K}^{(j+1)}\boldsymbol{u}^{(j+1)} + \boldsymbol{f}_{nl} 
\end{align}
Considering that the motion of the $j^{th}$ unit is decomposed over $M$ harmonics with frequencies $\omega, 2\omega, ... , M\omega$ and corresponding wave numbers $k, 2k, ... , Mk$:
\begin{align}
\label{eqn:harmonics}
&\boldsymbol{u}^{(j)}=A \left( \sum_{i=1}^{M}  \begin{bmatrix} c_{1i} \\ c_{2i} \end{bmatrix} \cos(i(jk L-\tau))+\begin{bmatrix} s_{1i} \\ s_{2i} \end{bmatrix} \sin(i(jk L-\tau))\right)
\end{align}
where $A$ is the wave amplitude and $\{c_{1i},s_{1i},c_{2i}\,s_{2i}\}$ are the normalized amplitude vectors of the $i^{th}$ unit for the main and internal resonators respectively. By changing the summation to matrix/vector multiplication, the harmonic motion defined in Eq. \ref{eqn:harmonics}, can be rewritten as:

\begin{align}
\label{eqn:disp}
&\boldsymbol{u}^{(j)}=A\begin{bmatrix} \boldsymbol{\psi}^{(j)}(\tau) & 0\\ 0 & \boldsymbol{\psi}^{(j)}(\tau) \end{bmatrix}_{2\times 4M}  
\begin{bmatrix} c_{11}, s_{11},c_{12}, s_{12}, \dots, c_{2M}, s_{2M}\end{bmatrix}^{-1}_{4M\times 1}
\end{align}
where \\
\begin{align}
\label{eqn:block}
\fontsize{10.4}{12}
\boldsymbol{\psi}^{(j)}(\tau) = \begin{bmatrix} cos(jk L-\tau)&sin(jk L-\tau)&...&cos(M(jkL-\tau))&sin(M(jkL-\tau))\end{bmatrix}_{2M\times 1}
\end{align}
In these equations, $j$ is the reference unit from which the wave propagates through neighboring units. Without loss of generality, j is set to 0, and $\boldsymbol{u}^{(j)}$, $\boldsymbol{u}^{(j-1)}$ and $\boldsymbol{u}^{(j+1)}$ are shown with $\boldsymbol{u}^{(0)}$, $\boldsymbol{u}^{(-1)}$ and $\boldsymbol{u}^{(1)}$. Eq. \ref{eqn:eqn_main3} can be rewritten using Eq. \ref{eqn:disp}.
\begin{align}
\label{eqn:eqn_main4}
&\omega^2 A \boldsymbol{M}\frac{d^2 \Psi^{(0)}}{d\tau^2} \boldsymbol{\lambda}+ \boldsymbol{F_{int}}(\boldsymbol{u}^{(0)},\boldsymbol{u}^{(-1)},\boldsymbol{u}^{(+1)}) = 0
\end{align}
Where,
\begin{align}
\Psi^{(0)}=\begin{bmatrix} \boldsymbol{\psi}^{(0)}(\tau) & 0\\ 0 & \boldsymbol{\psi}^{(0)}(\tau) \end{bmatrix}
\end{align}
\begin{align}
\boldsymbol{\lambda}=\begin{bmatrix} c_{11} &s_{11}&...&c_{1M}&s_{1M}&c_{21}&s_{21}&...&c_{2M}&s_{2M}\end{bmatrix}^T
\end{align}
The maximum element of $\boldsymbol{\lambda}$ is normalized to unity.
The differential equation (\ref{eqn:eqn_main4}) is converted to a set of $4Q$ nonlinear algebraic equations using Galerkin's projection \cite{bobylev2012approximation}:
\begin{align}
\label{eqn:eqn_main5}
&\int_{0}^{2\pi}   \Psi^{{(0)}^T} (\omega^2 A \boldsymbol{M}\frac{d^2 \Psi^{(0)}}{d\tau^2}\boldsymbol{\lambda}+ \boldsymbol{F_{int}}(\boldsymbol{u}^{(0)},\boldsymbol{u}^{(-1)},\boldsymbol{u}^{(+1)}))d\tau = 0
\end{align}
There are a total of $4M+2$ variables: $\boldsymbol{\lambda}$,  wave amplitude $A$, wave number $k$ and frequency $\omega$. In order to find the dispersion relation between the wave number and frequency, the amplitude  and the wave number are fixed  thus leading to $4M$ unknowns. We can now solve for  $\omega$ and $\boldsymbol{\lambda}$ for a given amplitude and wave number. 

\begin{figure}
	\centering
	\begin{subfigure}[b]{0.31\textwidth}
		\includegraphics[width=\linewidth,keepaspectratio]{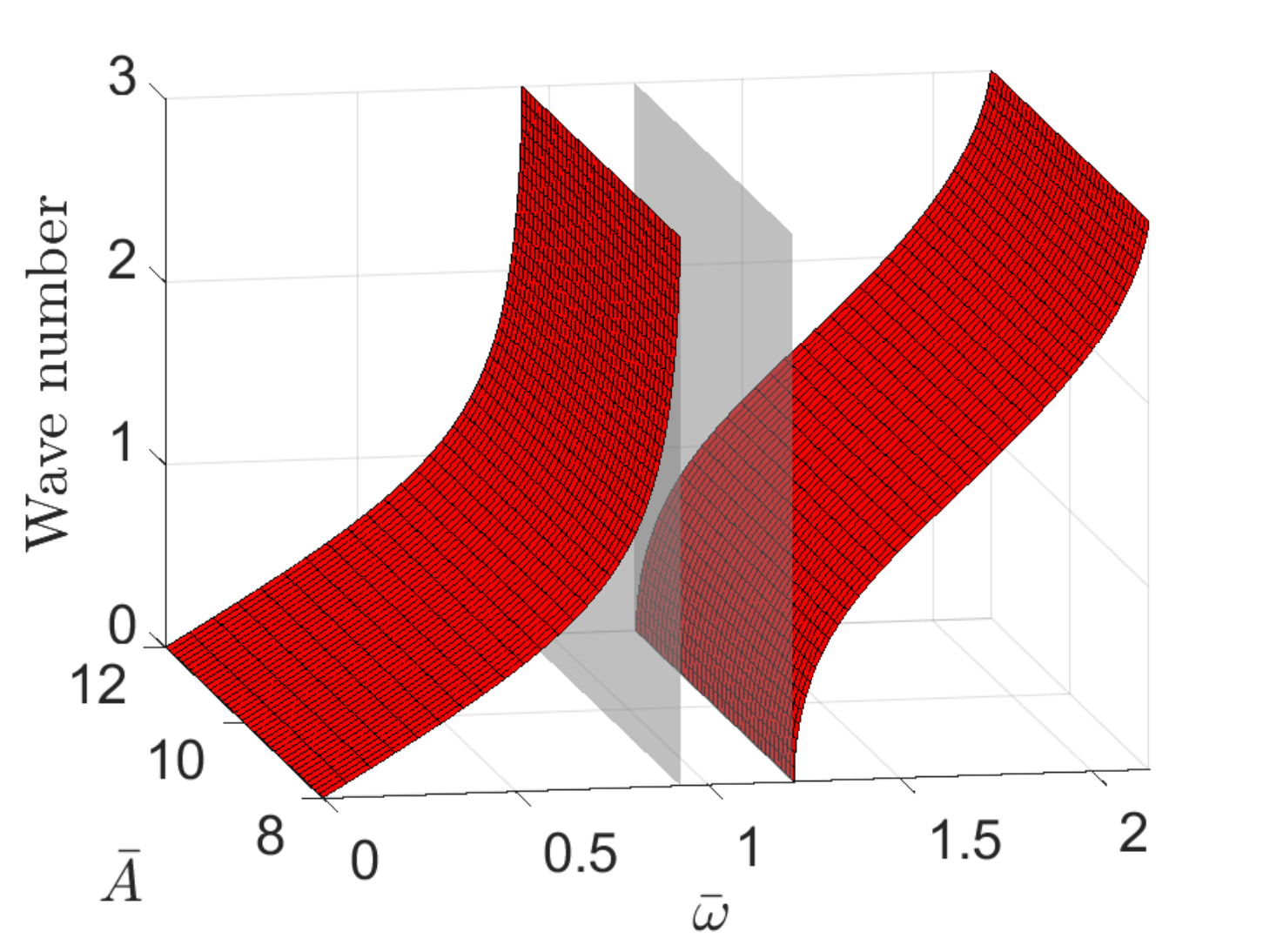}
		\caption{}
		\label{fig:bandgap_3D_a}
	\end{subfigure}
	\begin{subfigure}[b]{0.31\textwidth}
		\includegraphics[width=\linewidth,keepaspectratio]{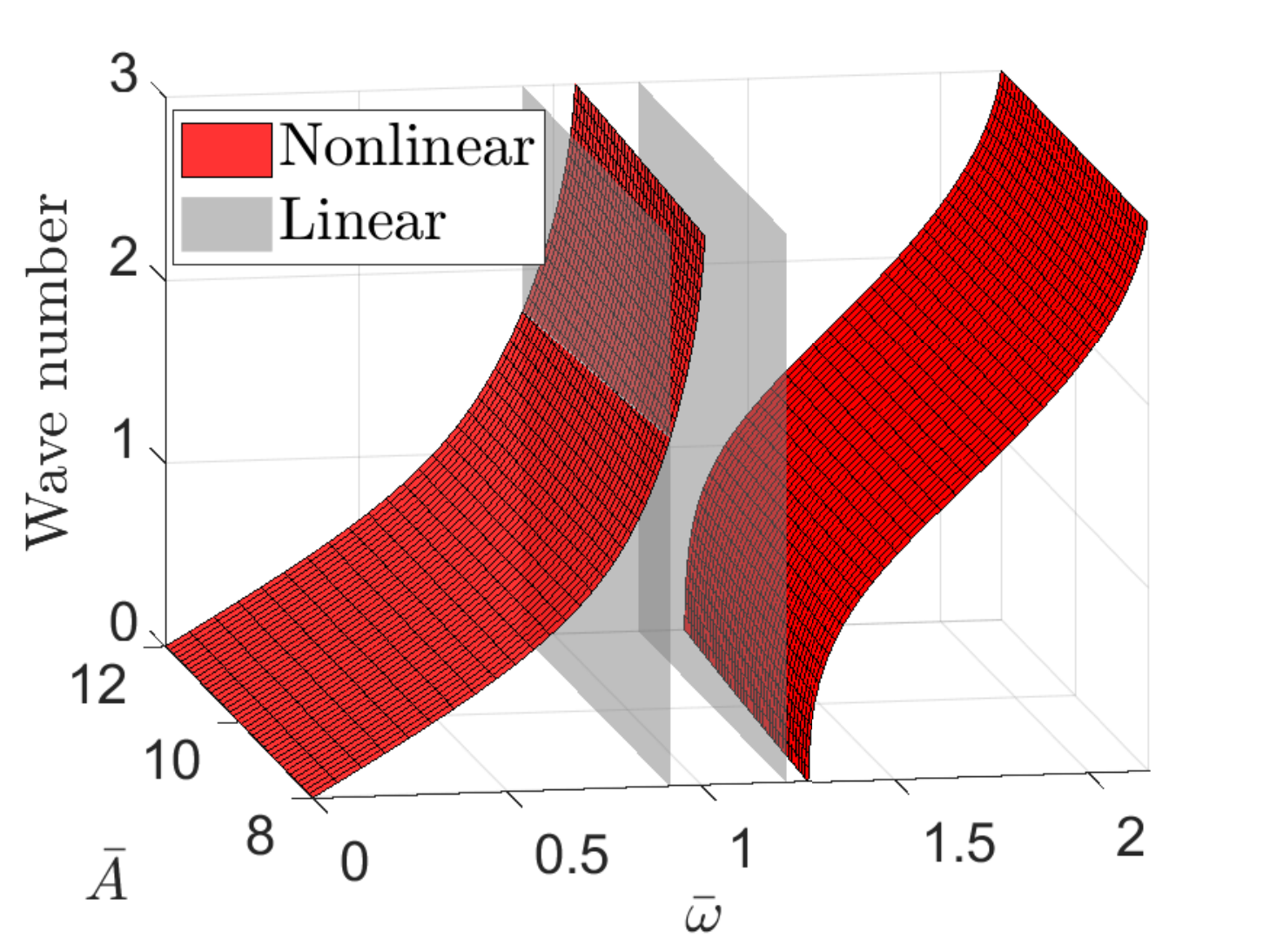}
		\caption{}
		\label{fig:bandgap_3D_b}
	\end{subfigure}
	\begin{subfigure}[b]{0.31\textwidth}
		\includegraphics[width=\linewidth,keepaspectratio]{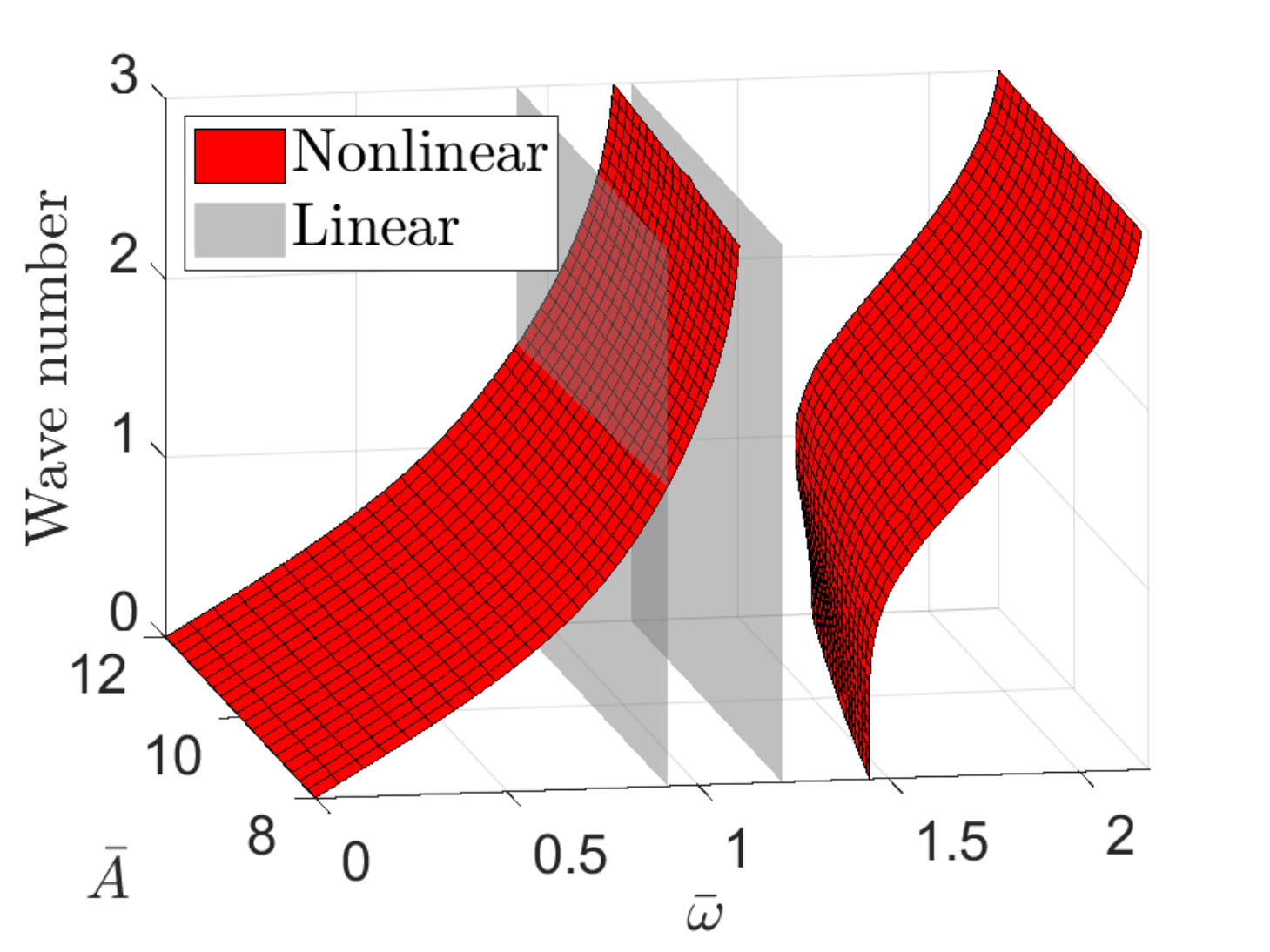}
		\caption{}
		\label{fig:bandgap_3D_c}
	\end{subfigure}
	\caption{\small{Dispersion behavior of a chain of resonators as a function of amplitude and nonlinear stiffness for $\alpha=0.5$ and $\epsilon=0.5$. (a) $\eta=0$. (b): $\eta=1\times 10^{-4}$. (c): $\eta=5\times 10^{-4}$.}\label{fig:bandgap_3D}}
\end{figure}

Figure \ref{fig:bandgap_3D} provides an example of dispersion with the corresponding band gaps. It shows how the dispersion and the band gaps change as the function of the amplitude and the nonlinear stiffness of the internal resonator $k_n$.

\subsection{Influence of nonlinear stiffness. Discontinuous behavior.}
\begin{figure}[!h]
	\centering
	\includegraphics[width=0.48\linewidth]{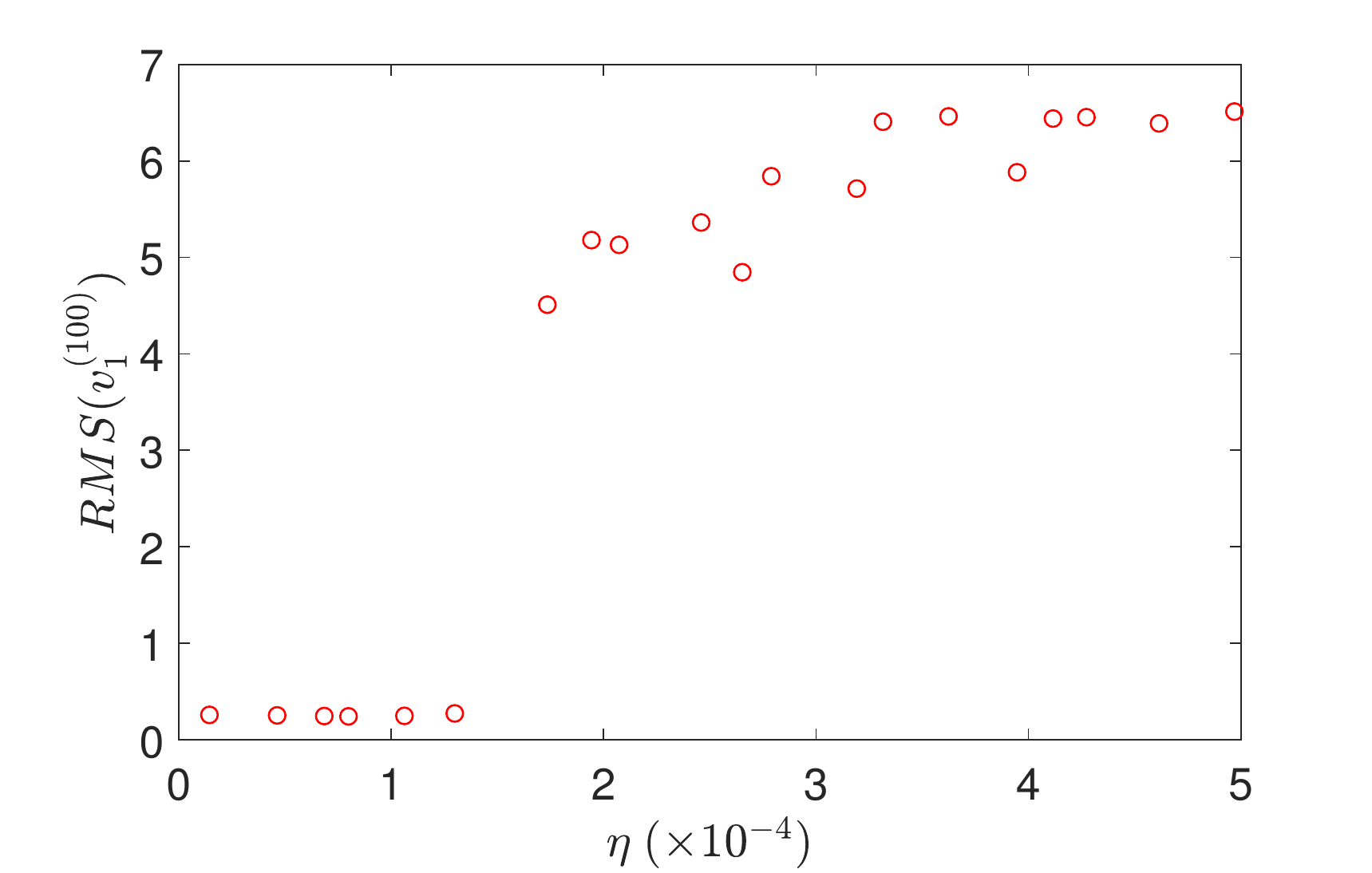}
	\caption{\small{RMS response of the $100^{th}$ unit  of a 1000 unit chain subjected to a harmonic excitation with respect to  the nonlinear stiffness $\eta$ of the internal resonators.  $\alpha=0.5$, $\epsilon=0.5$, $\xi_1 = 0$, $\xi_2 = 0$, $\bar{A}=10$ and $\bar{\omega}=1$.}}
	\label{fig:disc}
\end{figure}

As shown in the previous section, the dispersion relation, and therefore the dynamic behavior of the chain, is dependent on the amplitude and also on the nonlinear stiffness of the internal resonator. Another important feature of the dynamic response is that it can exhibit discontinuities. {Figure \ref{fig:disc} depicts an example of Root Mean Square (RMS) response as a function of the nonlinear stiffness of the internal resonator for a mass-in-a-mass design.} The response exhibits a discontinuity. Such discontinuity will reduce as the damping increases. It is noteworthy that when the internal resonator has non-zero damping $\xi_2 \neq 0$, it might behave like a nonlinear energy sink (NES) which has shown to exhibit an activation threshold  \cite{boroson2017optimization, boroson2017stochastic} leading to a discontinuous behavior.  \\

The discontinuities make the chain potentially highly sensitive to uncertainty. In addition, they present a challenge from the optimization point of view as the approximation of the response becomes tedious. The need to include uncertainties and tackle discontinuities has led to the development of specifically tailored stochastic optimization methods that were used in the context of the optimization of NESs \cite{boroson2017stochastic}. The optimization formulation and a summarized version of the algorithm are provided in the following section.\\

\section{Optimization formulation and algorithm}
\label{sec:optim_formul}
\subsection{Optimization formulation and calculation of expected value}
This article aims at optimizing a chain of resonators through the following stochastic optimization formulation:

\begin{align}
\label{eqn:stoch_optim}
\min_{\boldsymbol{\mu}_d}
& \quad  \mathbb{E} (F(\mathbf{X}_d,\mathbf{X}_a))\\
s.t.&\quad \boldsymbol{\mu}_d^{min}\leq\boldsymbol{\mu}_d\leq \boldsymbol{\mu}_d^{max}\notag
\end{align}

where, $\mathbb{E}$ is the expected value. $F$ is a quantity to optimize such as an RMS displacement at a given point on the chain.  $\mathbf{X}_d$ and $\mathbf{X}_a$ are respectively the (random) design and ``aleatory" variables. The aleatory variables are defined as variables which are not controllable. $\boldsymbol{\mu}_d$ are the means of the design variables and $\boldsymbol{\mu}_d^{min}$ and $\boldsymbol{\mu}_d^{max}$ are respectively their optimization lower and upper bounds.  The linear and nonlinear stiffnesses of the resonators are examples of $\mathbf{X}_d$, and the amplitude and frequency of excitation are examples of $\mathbf{X}_a$. 


\subsection{Field description}
\label{sec:field}
Chains of resonators can include a very large number of units in which the properties such as stiffness and mass should be adjusted. Therefore,  it would be impossible to optimize units individually. In this paper, we propose to reduce the dimensionality of the problem by describing the unit features (e.g., nonlinear stiffness) through a field controlled with a handful of coefficients, which become the new optimization variables of the problem. This approach makes the problem scalable by significantly reducing the number of optimization variables. The proposed field $\theta$ is expanded over $m$ ``basis'' functions:
\begin{equation}
\theta(x) = a_1 + \sum_{j=1}^{m} b_{j} \phi (c_{j} x) \quad
\label{eqn:field_general}
\end{equation}
where $x$ represents the position of a unit on the chain, and $a_1$, $b_j$ and $c_j$ are the coefficients, which will become optimization variables. {The rationale for the proposed field form is that, from a design standpoint,  it enables the representation of a range of periodic as well non-periodic structures in the case of a finite chain.} If only one term is used and $\phi(x)=sin^2(x)$, as used in this work, the field becomes:
\begin{equation}
\theta(x) = a_1 +  a_{2} sin^2(a_{3} x ) 
\label{eqn:field}
\end{equation}
where, for ease of notation, $a_{2}\equiv b_{1}$ and  $a_{3}\equiv c_{1}$. If one term is used, defining a property of the chain using the field enables one to reduce the number of optimization variables to three. {Note that the purpose of squaring the basis function is to obtain non-negative values from the field description.}

\begin{figure}
	\centering
	\begin{subfigure}[b]{0.49\textwidth}
		\includegraphics[width=\linewidth,keepaspectratio]{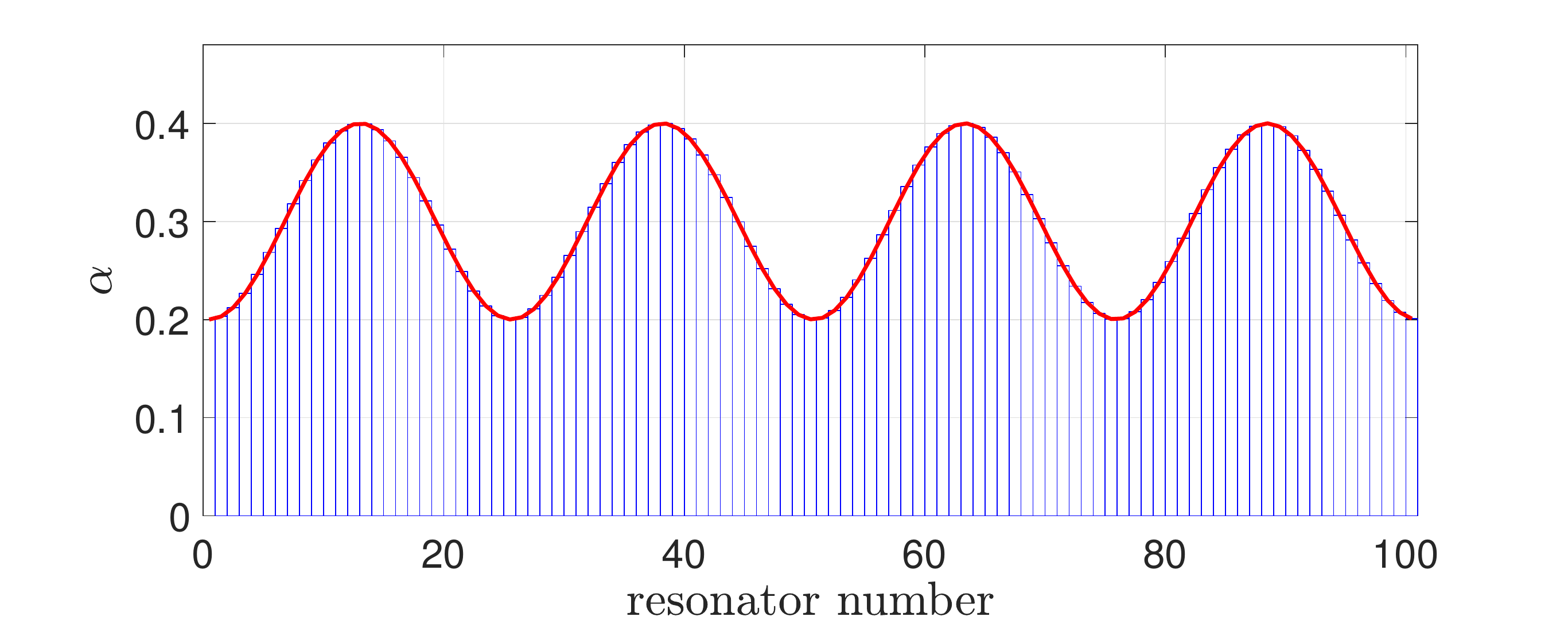}
		\caption{}
		\label{fig:field_a}
	\end{subfigure}
	\begin{subfigure}[b]{0.49\textwidth}
		\includegraphics[width=\linewidth,keepaspectratio]{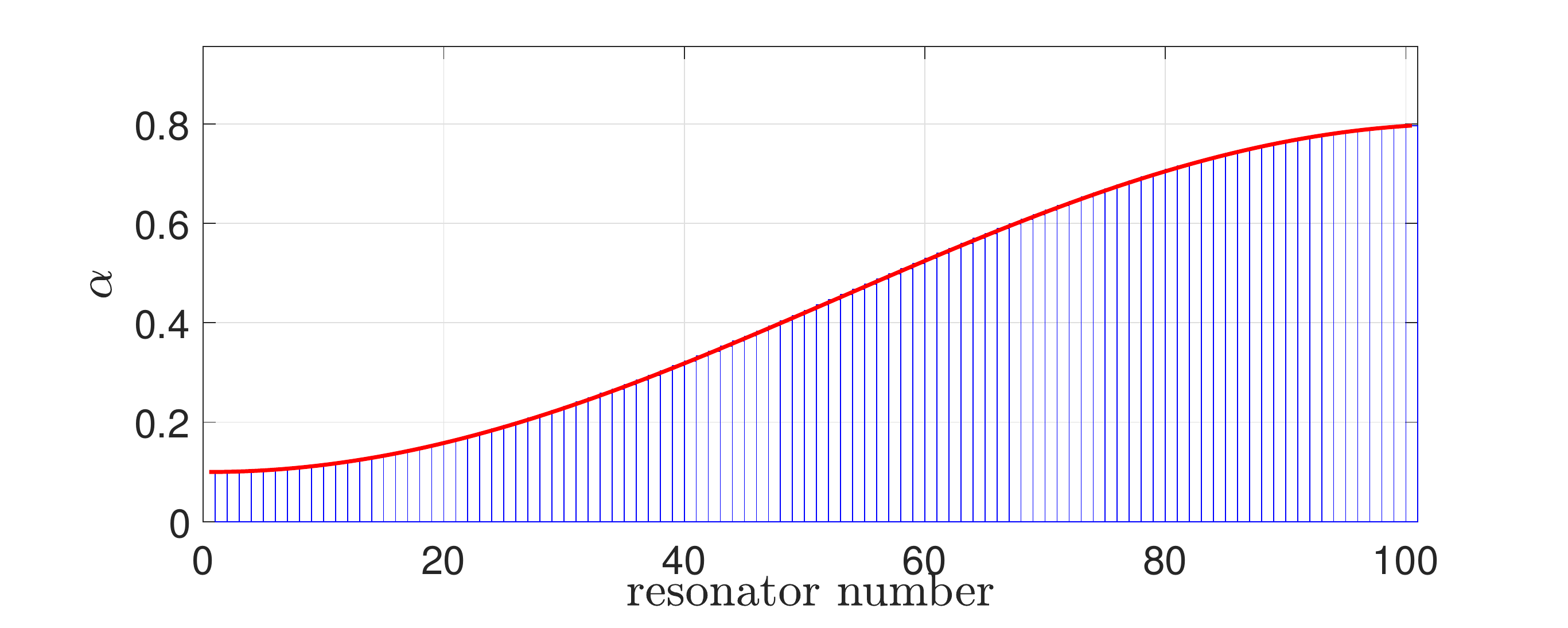}
		\caption{}
		\label{fig:field_b}
	\end{subfigure}
	\caption{\small{Variation of the stiffness ratio $\alpha$ through a chain of 100 units defined using the field function. (a): the coefficients are $a_1=0.2$, $a_2=0.2$, $a_3=0.025$ and $\alpha$ varies periodically. (b): the coefficients are $a_1=0.1$, $a_2=0.7$, $a_3=0.003$ and $\alpha$ increases monotonically.}}\label{fig:field}
\end{figure}

For two specific choices of coefficients $a_1$, $a_2$ and $a_3$, the variations of the stiffness ratio $\alpha$ in a chain of 100 units are illustrated in Figure \ref{fig:field}. For the first choice of coefficients ($a_1=0.2$, $a_2=0.2$, $a_3=0.025$), the variation of $\alpha$ has become periodic through the chain. However, for the second set of coefficients ($a_1=0.1$, $a_2=0.7$, $a_3=0.003$), $\alpha$ increases monotonically.

\section{Optimization algorithm}
The stochastic optimization problem in Eq. \ref{eqn:stoch_optim} is solved using a specific optimization algorithm \cite{boroson2017stochastic} whose key elements and steps are described below.

\begin{figure}[!h]
	\centering
	\includegraphics[width=0.65\linewidth]{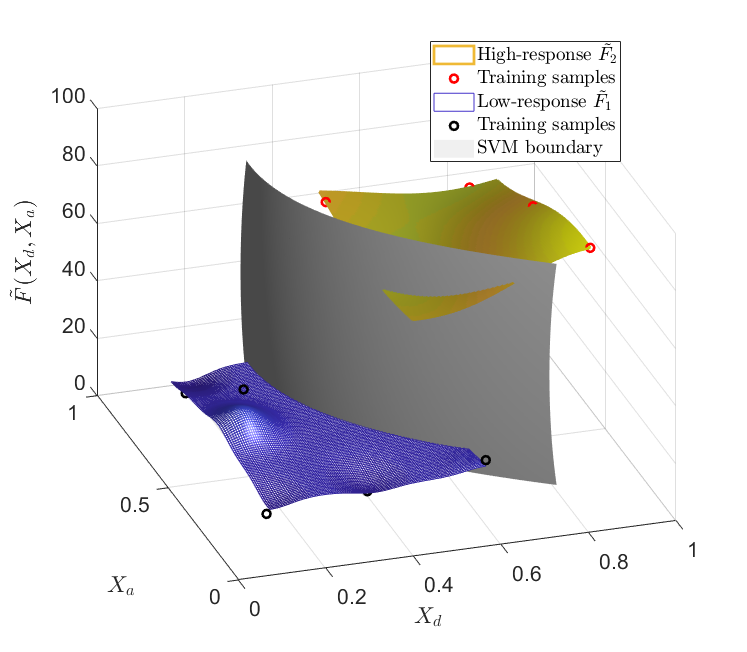}
	\caption{\small{An example of SVM boundary with two Kriging approximations trained using two clusters.}}
	\label{fig:kriging_svm}
\end{figure}

In order to optimize the problem, a surrogate approximation of the response must be used due to the computational burden of each simulation. However, the use of a single meta-model such as Kriging constructed with a design of experiments would not lead to a good approximation because of the presence of potentially large discontinuities. {For this reason, a clustering technique such as K-Means \cite{Hartigan1979} is used to detect the groups of responses that are distinct because of the response discontinuity. The samples corresponding to the two clusters of responses are then labeled into classes and are used to train an SVM classifier \cite{christianini2000support}. The SVM provides an explicit boundary that splits the space defined by random design and aleatory variables into two regions. The response in each region is assumed continuous and can therefore be approximated using a separate Kriging approximation.} An example of two Krigings, trained using two clusters, and an SVM boundary, which splits this space into two regions with distinct levels of response is illustrated in Figure \ref{fig:kriging_svm}.

Sub-problems based on two Krigings and an SVM are refined iteratively through a specific adaptive sampling scheme which accounts for the random variable distributions \cite{lacaze2014generalized}. Monte-Carlo simulations over the Krigings are used to approximate the expected values. In addition, the method makes use of the availability of the Kriging prediction variance to find the aleatory variables during the optimization. Note that a global and zero-order optimizer is needed -- the problem is non-differentiable at the SVM boundary-- to solve the problem. In this work, Particle Swarm Optimization (PSO) is used  \cite{kennedy2011particle}. The use of PSO does not involve large computational time since all PSO function evaluations are based on the Kriging approximation. The following describes the main steps of the algorithm.

\subsection{Clustering}
The first step of this optimization algorithm is to detect the discontinuity in the response and split the samples into two classes with high and low responses. In this study, the K-Means algorithm is used to cluster the samples based on the RMS displacement. 

\subsection{Support vector machine (SVM)}
Using the two clusters of samples, the regions corresponding to high and low RMS values can be identified in the space of parameters. The idea is to construct a boundary that splits the parameters space into two regions which are associated with low and high RMS values. For this purpose, an SVM classifier is used which is trained based on the labeled samples and provides an explicit expression of the boundary in terms of the parameters. Given a set of $N$ training samples $\mathbf{x}_{i}$ in an $n$-dimensional space and the corresponding class labels $y_i=\pm 1$, which in our case is provided by the K-means clustering, an SVM boundary is given as:

\begin{equation}
\label{decisionfinal}
s(\mathbf{x})=b+\sum_{i=1}^{N}\lambda_{i}y_{i}K(\bf{x}_{i},\bf{x})=0
\end{equation}

\noindent where $b$ is a scalar referred to as the bias, $\lambda_{i}$ are Lagrange multipliers obtained from the quadratic programming optimization problem used to construct the SVM, and $K$ is a kernel function. The classification of any arbitrary point $\mathbf{x}$ is given by the sign of $s(\mathbf{x})$.  The kernel function $K$ can have several forms, such as polynomial or Gaussian radial basis kernel, which is used in this article: \begin{equation}
\label{eqgaussian}
K(\mathbf{x}_{i},\mathbf{x}_{j})=\exp{ \left( - \frac{\left|\left|\mathbf{x}_i-\mathbf{x}_j \right|\right|^2 }{2\sigma^2}            \right)}
\end{equation}
where $\sigma$ is the width parameter.

\subsection{Refinement of the SVM boundary using adaptive sampling}
Because the constructed SVM boundary might not be accurate around the optimization iterate, an adaptive sampling scheme is used to refine the boundary.  The sampling algorithm is described in detail by Lacaze and Missoum \cite{lacaze2014generalized}. A fundamental aspect of the algorithm is the selection of samples in the sparse regions of the space (i.e., as far away as possible from the existing samples) and also in the regions of highest probability of misclassification by the SVM. The latter criterion is obtained by locating the samples on the SVM. These samples are found by solving the following global optimization problem:

\begin{align}
\max_\mathbf{x} & \quad \frac{1}{d} \log\left(\textbf{f}_{\textbf{X}}(\textbf{x})\right) -\frac{1}{p} \log\left(\sum_{i=1}^{N_s}  ||\mathbf{x}-\mathbf{x}_i||^{-p}\right) \\
\text{s.t.}&\quad  s(\mathbf{x})=0 \notag
\end{align}

\noindent where $\textbf{f}_X$ is the d-dimensional standard normal joint probability density function, $N_s$ is the number of existing samples, and $p\gg1$. This approach selects points on the boundary of the SVM weighted by the known distributions of the design parameters. 

\subsection{Expanded space}
\begin{figure}[!h]
	\centering
	\includegraphics[width=0.60\linewidth]{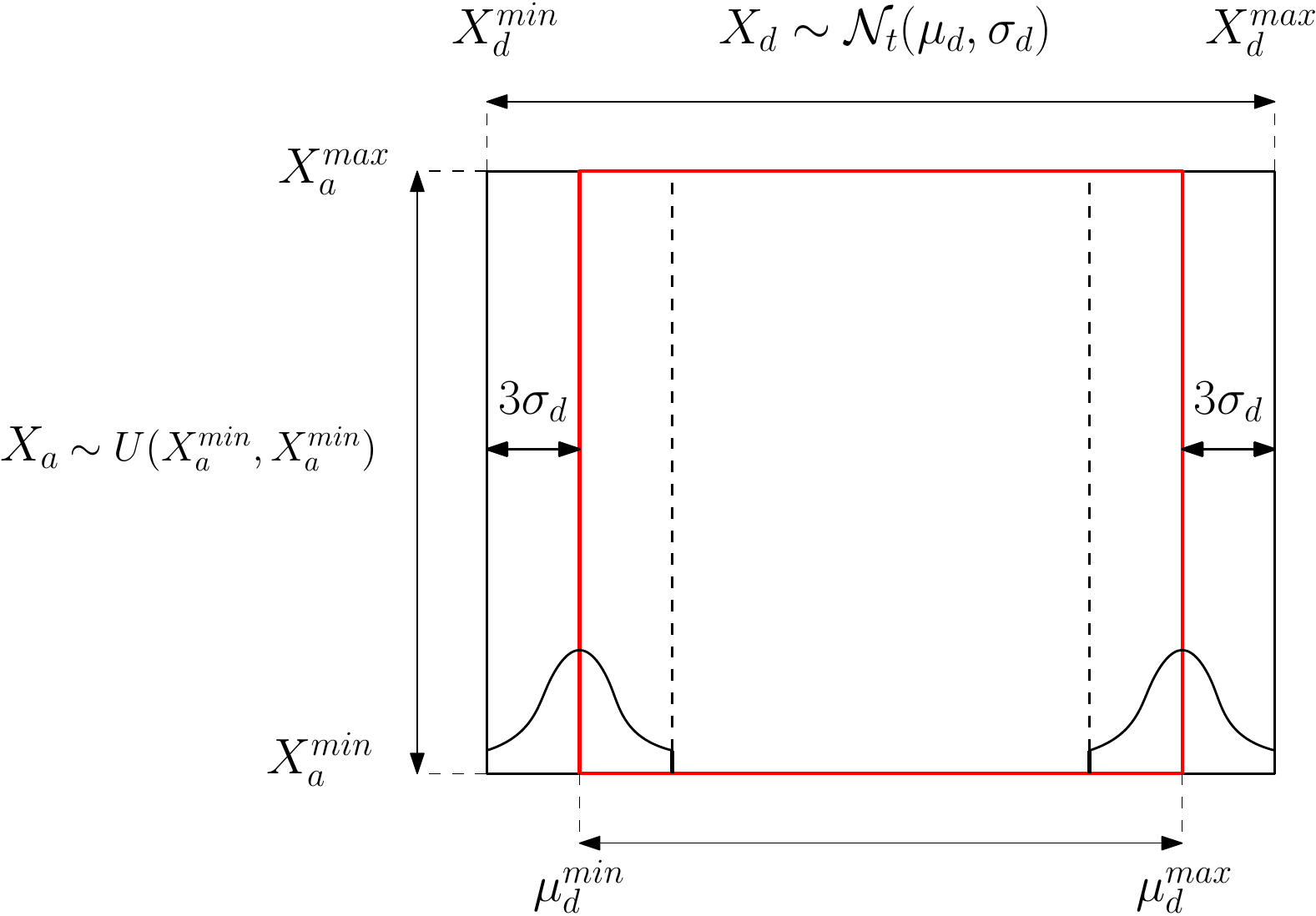}
	\caption{\small{Conceptual representation of the boundaries for the optimization and the expanded space where the surrogate models are defined.}}
	\label{fig:ex_space}
\end{figure}

As mentioned in Section \ref{sec:optim_formul}, the lower and upper bounds for the means of the design variables are $\boldsymbol{\mu}_d^{min}$ and $\boldsymbol{\mu}_d^{max}$ which specify the space where the optimization problem (Eq. \ref{eqn:stoch_optim}) is solved. However, because of the distributions defined for the design variables to perform the Monte-Carlo sampling, the samples are created in an expanded space whose boundaries for the design variables are $\mathbf{X}_d^{min}$ and $\mathbf{X}_d^{max}$ . Thus, the surrogate models are created and refined in this expanded space to ensure that the calculation of expected value based on the Monte-Carlo samples is accurate.
$\mathbf{X}_d^{min}$ and $\mathbf{X}_d^{max}$ are calculated respectively based on $\boldsymbol{\mu}_d^{min}$ and $\boldsymbol{\mu}_d^{max}$, and their corresponding distributions. For instance, in case of using a truncated normal distribution bounded by $\mu\pm 3\sigma$, $\mathbf{X}_d^{min}$ and $\mathbf{X}_d^{max}$ are respectively $\boldsymbol{\mu}_d^{min}-3\boldsymbol{\sigma}_d$ and $\boldsymbol{\mu}_d^{max}+3\boldsymbol{\sigma}_d$, as shown in Figure \ref{fig:ex_space}.

\subsection{Approximation of response using two Kriging approximations. Computation of the expected value}
In order to approximate the response over the entire space (including both design and aleatory variables), two separate Kriging approximations are constructed. Each one is built based on one of the clusters as identified through clustering. The optimization to minimize $\mathbb{E} (F)$ is performed using these two approximations. Monte-Carlo simulations are used to compute the expected value around the optimization iterate. In order to decide which Kriging approximation will be used to evaluate the Monte-Carlo samples, the SVM described above is used. Based on $N_{MC}$ Monte-Carlo samples $\mathbf{x}_i$, the approximated expected value can therefore be written as:

\begin{align}
&\label{eq:MCprob}
\mathbb{E} ({F})\approx\frac{1}{N_{MC}}\sum_{i=1}^{N_{MC}}  I(\mathbf{x}_i\in\Omega)\tilde{F}_1(\mathbf{x}_i) + I(\mathbf{x}_i\in\bar\Omega)\tilde{F}_2(\mathbf{x}_i)\\
&\Omega=\{\mathbf{x}|s(\mathbf{x})<0\} \notag
\end{align}
where $I$ is an indicator function defined as:
\[ I_s(\mathbf{x}\in\Omega)= \left\{
\begin{array}{l l}
1 & \quad \text{if $s(\mathbf{x}) \leq 0$}\\
0 & \quad \text{if $s(\mathbf{x}) > 0$}
\end{array} \right.\]
$\tilde{F}_1$ and $\tilde{F}_2$ are the two Kriging approximations of response corresponding to the two regions with low and high values. {Note that the Monte-Carlo simulations based on Kriging approximations are computationally efficient}.

\subsection{Sampling of the aleatory space. Update of SVM and Kriging}
\label{sec:maxvar}
The two Kriging approximations and the SVM are updated iteratively using the generalized max-min sample as well as the optimization iterate. However, the optimization is only performed in the space of design variables whereas the SVM and the two Krigings are constructed in the design and aleatory space. In order to provide the optimization iterate with an aleatory component, the readily available prediction variance of the Kriging approximations is used. Specifically, the aleatory variables at the optimization iterate are chosen as the maximizer of the Kriging prediction variance, which can be considered where the Kriging is the least accurate. For a given design iterate $\boldsymbol{\mu}_d^{(k)}$,  the aleatory variables can be found through the following optimization problem:
\begin{align}
\label{eq:maxvar}
\mathbf{x}_a^{(k)}=\underset{\mathbf{x}_a}{\arg\max} & \quad \sigma_{K}(\boldsymbol{\mu}_d^{(k)},\mathbf{x}_a)\\\notag
\text{s.t.} & \quad \mathbf{X}_a^{min} \leq \mathbf{x}_a \leq \mathbf{X}_a^{max}\notag
\end{align}
where $\sigma_{K}$ is the predicted variance of the Kriging approximations at the $k^{th}$ iteration. This sample will be referred to as the ``maximum variance sample". Because the algorithm uses two Kriging approximations, the search is performed over the two surrogates. 

\section{Results}
\label{sec:results}
The stochastic optimization of a chain composed of 1000 units is considered. {For the results presented in this section, the chain is subjected to forced harmonic displacement excitation(s) applied to the main resonator of unit 1.} The response considered ($F$ in Equation \ref{eqn:stoch_optim}) is the RMS displacement of unit 100. In each optimization, the random variables are:
\begin{itemize}
	\item The design variables $\mathbf{X}_d$. In all problems, truncated normal distributions $\mathbf{X}_d \sim \mathcal{N}_t(\boldsymbol{\mu}_d,\boldsymbol{\sigma}^2)$ are used. The bounds of the distributions are $\boldsymbol{\mu}_d \pm 3\boldsymbol{\sigma}_d$. The standard deviation is calculated based on a coefficient of variation of 2\%. 
	\item The amplitude $\bar{A}$  and frequency $\bar{\omega}$ of excitation are chosen as aleatory variables. They follow uniform distributions with bounds of $\mathbf{X}_a^{min}$ and $\mathbf{X}_a^{max}$. 
\end{itemize}

\noindent Three optimization problems are presented:
\begin{itemize}
	\item {\em Optimization problem 1}\\
	Random design variable: nonlinear stiffness ${X}_d=\{\eta\}$ constant over the chain. One excitation (frequency and amplitude are aleatory variables). Total of three variables.
	\item {\em Optimization problem 2}\\
	Random design variable: stiffness ratio, mass ratio and nonlinear stiffness, $\mathbf{X}_d=\{\alpha, \epsilon, \eta\}$ constant over the chain. Two excitations (frequencies and amplitude are aleatory variables). Total of 6 variables.
	\item {\em Optimization problem 3}\\
	Same as optimization problem 2 with design variables represented using a field. Total of 12 variables.
\end{itemize}

\noindent A linear undamped chain is used as a reference design (Table \ref{tab:reference}).
\begin{table}[!h]
	\centering
	\caption{\small{Reference design: linear undamped chain.}}
	\label{tab:reference}
	\begin{tabular}{llllll}
		\Xhline{2pt}
		\centering
		$\alpha$&$\epsilon$ &$\eta$  &$\xi_1$ &$\xi_2$  \\
		\hline
		\centering
		0.5&0.5& 0 & 0 & 0    \\
		\Xhline{2pt}
	\end{tabular}
\end{table}

\noindent For the interested reader, the appendix provides details related to the accuracy of the surrogate-based optimization and the variance of the results due to Monte-Carlo sampling related to the three optimization problems. 
\subsection{Optimization problem 1. 1 design variable. 2 aleatory variables}
\begin{figure}[!h]
	\centering
	\includegraphics[width=0.50\linewidth]{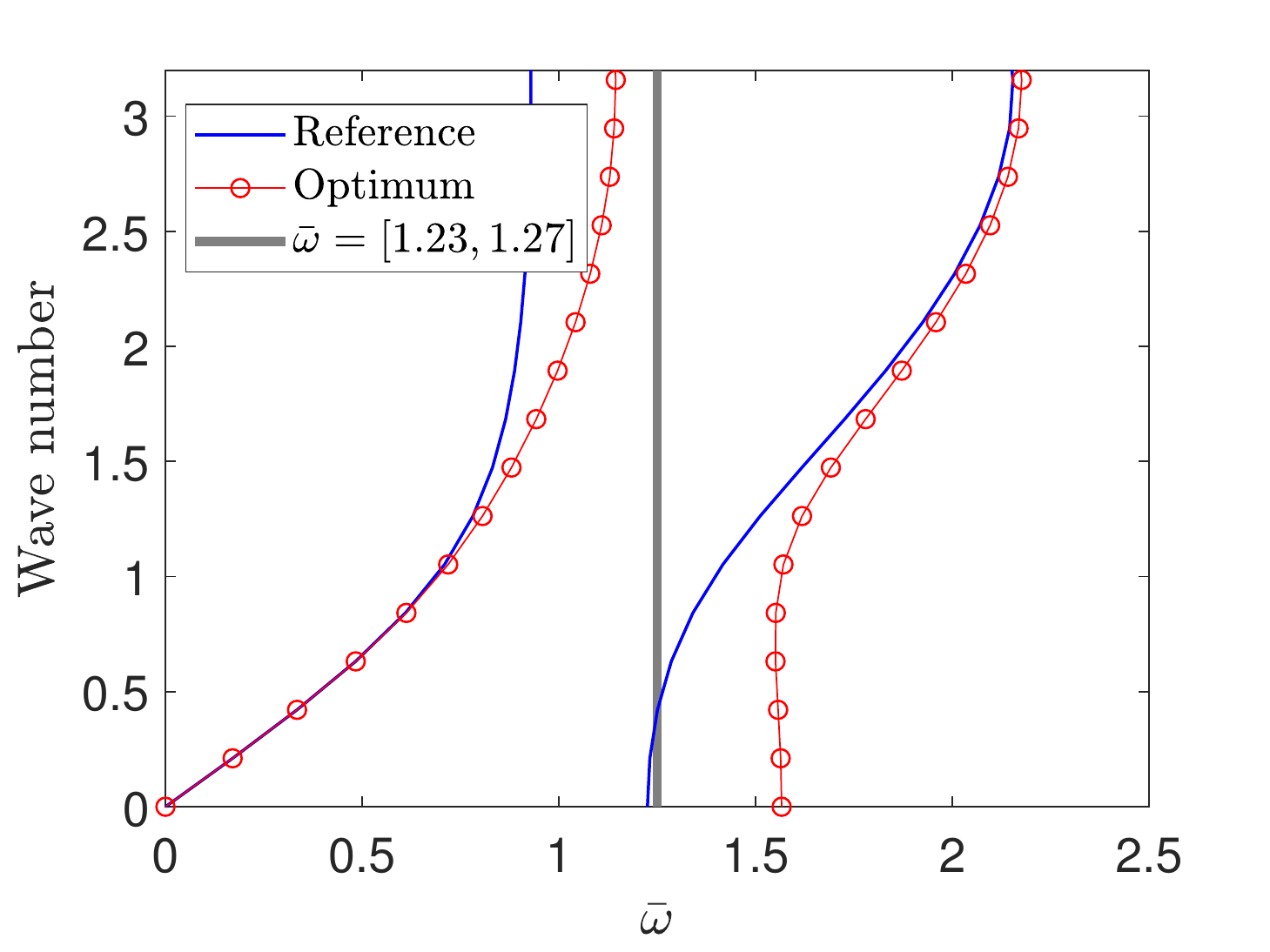}
	\caption{\small{Problem 1. Comparison of the dispersion diagram of the reference design and the optimal design {($\eta = 4.72 \times 10^{-4}$)} for amplitude of $\bar{A}=10$.}}
	\label{fig:optim1_HB}
\end{figure}

In this first problem, the random design variable, non-dimensional nonlinear stiffness ${X}_d=\{\eta\}$, is assumed constant over the chain. The excitation amplitude and frequency are considered as aleatory variables. Their ranges are $\bar{A}=[8,12]$ and $\bar{\omega}=[1.23,1.27]$ respectively. As shown in Section \ref{sec:HM}, changes in the nonlinear stiffness and the amplitude of excitation have an influence on the band gap. In this problem, the frequency range falls outside the band gap of the reference design as depicted in Figure \ref{fig:optim1_HB}.  The optimization problem reads:
\begin{align}
\label{eqn:opt1}
\min_{\mu_{\eta}}
& \quad  \mathbb{E} ({F}(\eta, \bar{A},\bar{\omega}))\\ \notag
s.t.
&\quad 10^{-6} + 3 \sigma_{\eta} \leq \mu_{\eta}\leq 5 \times 10^{-4} - 3 \sigma_{\eta} \notag
\end{align}

The variable distributions are:
\begin{align}
F\equiv RMS(v_1^{(100)}) \notag \\
{X}_d \sim \mathcal{N}_t(\mu_{\eta},\sigma_{\eta}^2), \; \sigma_{\eta} = 0.02\mu_{\eta}   \notag\\
\bar{A} \sim U(8,	12) \notag \\
\bar{\omega} \sim U(1.23, 1.27) \notag
\end{align}

Results are listed in Table \ref{tab:opt1_result}. The table also provides the relative reduction in RMS displacement value compared to the reference configuration for which the RMS is 5.29. For this problem, the reduction in RMS displacement value is 82.8\%.

\begin{table}[!h]
	\small
	\centering
	\caption{\small{Problem 1. Optimization results.}}
	\label{tab:opt1_result}
	\begin{tabular}{ccc}
		\Xhline{2pt}
		\centering
		$\mathbf{\mu}_{\eta}$&$\mathbb{E} ({F})$ & Rel. red. ref.   \\
		\hline
		\centering
		$4.72 \times 10^{-4}$&0.910& 82.8\%  \\
		\Xhline{2pt}
	\end{tabular}
\end{table}

The dispersion diagram of the optimized chain is provided in Figure \ref{fig:optim1_HB}. For comparison, the figure also provides the dispersion behavior for the reference (linear, undamped) design for the mean value of amplitude. It can be seen that the small excitation frequency range, which initially fell outside the band gap is now located within the band gap of the optimal nonlinear chain. It can be observed that the effect of the nonlinearity is to shift the band gap to higher frequencies to mitigate the vibration.

\begin{figure}
	\centering
	\begin{subfigure}[b]{0.4\textwidth}
		\includegraphics[width=\linewidth,keepaspectratio]{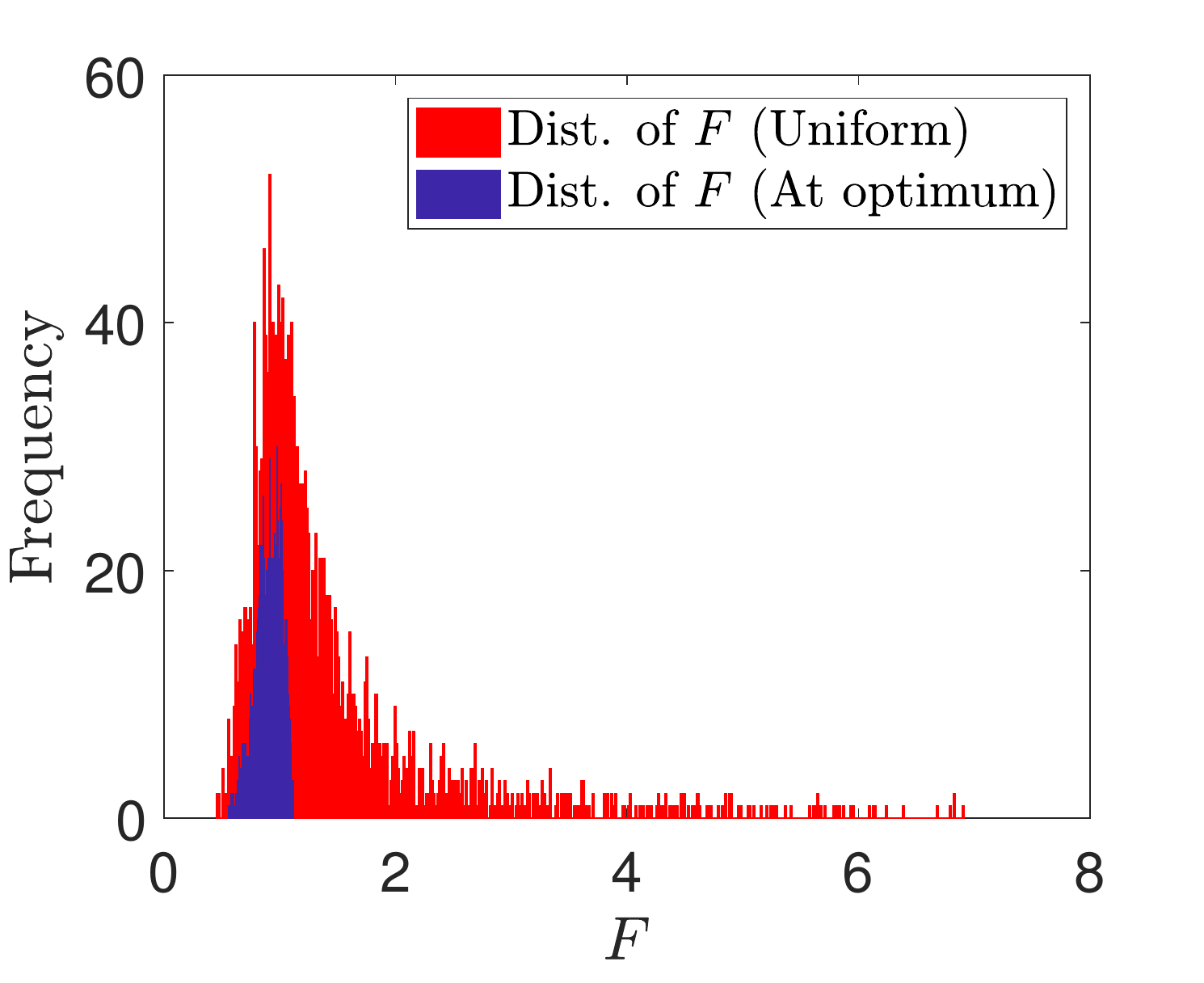}
		\caption{}
		\label{fig:hist_unif_problem1_a}
	\end{subfigure}
	\begin{subfigure}[b]{0.4\textwidth}
		\includegraphics[width=\linewidth,keepaspectratio]{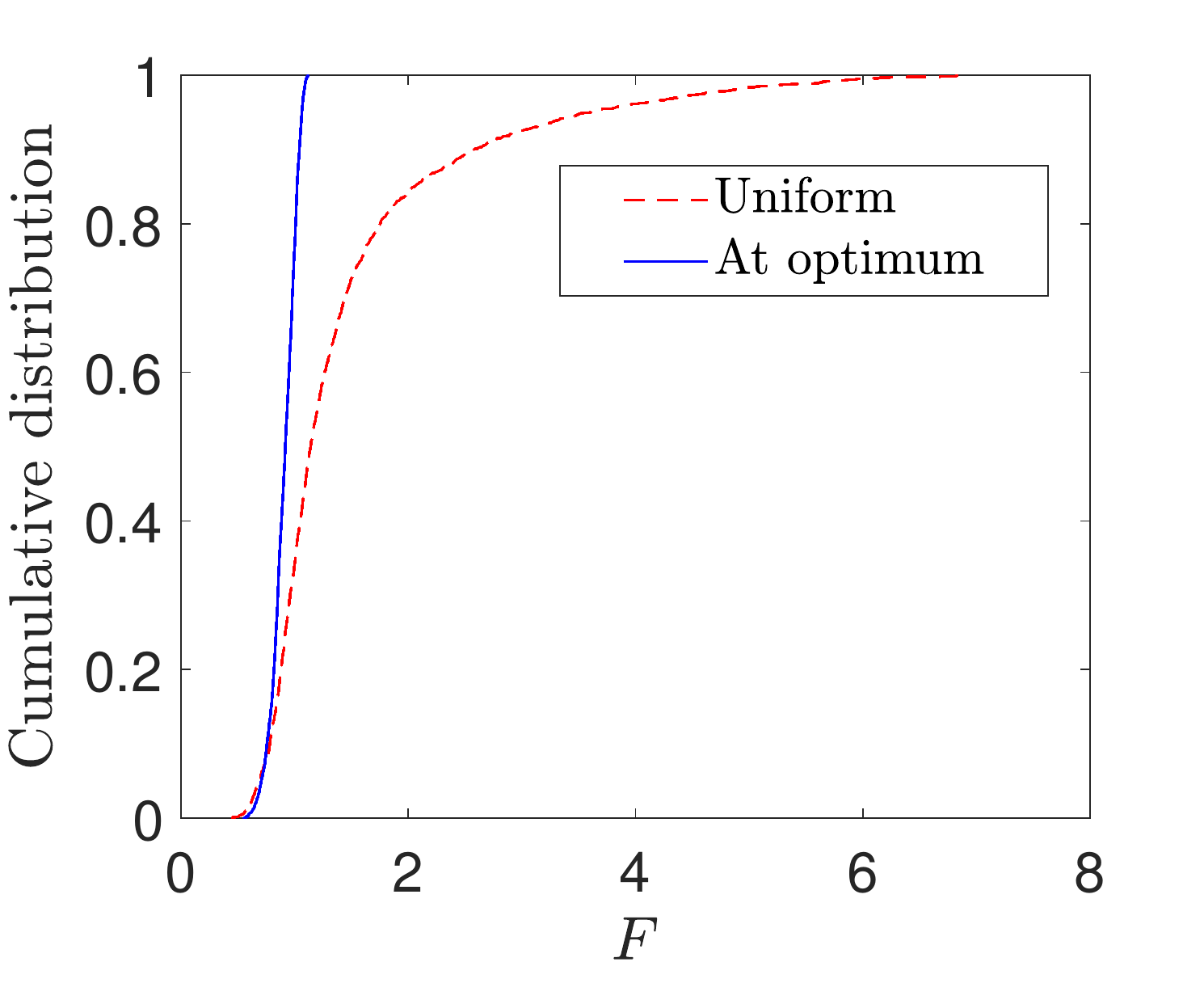}
		\caption{}
		\label{fig:hist_unif_problem1_b}
	\end{subfigure}
	\caption{\small{Problem 1. Comparing the distributions of RMS responses ($F$) around optimum vs. variables uniformly distributed. (a): Histograms of distributions. (b): Cumulative distributions.}}
	\label{fig:hist_unif_problem1}
\end{figure}

To further investigate the optimal design, the distribution of RMS response around the optimum is compared with the distribution of responses for any possible non-optimal design. {In order to evenly account for all possibilities of non-optimal designs, a uniform distribution is considered for all variables. As illustrated in Figure \ref{fig:hist_unif_problem1}(a), the variance of the response for the optimal design is significantly smaller than the non-optimal cases. This conclusion can also be reached by inspecting the cumulative distributions provided in Figure \ref{fig:hist_unif_problem1}(b). The results demonstrate that the ``robustness" of the optimal design and show that the non-optimal cases have a much larger spread toward the high RMS values.}

\subsection{Optimization problem 2. 3 design variables. 3 aleatory variables}
In the second problem, the optimization variables include stiffness ratio, mass ratio and nonlinear stiffness corresponding to the internal resonators, design variables $\mathbf{X}_d=\{\alpha, \epsilon, \eta\}$. The chain is subjected to a forced displacement excitation with two frequency ranges of $\bar{\omega}_1=[0.48,0.52]$ and $\bar{\omega}_2=[1.23,1.27]$. The optimization problem reads:

\begin{align}
\label{eqn:opt2}
\min_{\mu_{\alpha},\mu_{\epsilon},\mu_{\eta}}
& \quad  \mathbb{E} ({F}(\alpha,\epsilon,\eta, \bar{A},\bar{\omega}_1,\bar{\omega}_2))\\ \notag
s.t.
&\quad 10^{-6} + 3 \sigma_{\alpha} \leq \mu_{\alpha}\leq 0.5  - 3 \sigma_{\alpha} \notag\\
&\quad 0.25 + 3 \sigma_{\epsilon} \leq \mu_{\epsilon}\leq 2.5 - 3 \sigma_{\epsilon} \notag\\
&\quad 10^{-6} + 3 \sigma_{\eta} \leq \mu_{\eta}\leq 5 \times 10^{-4} - 3 \sigma_{\eta} \notag
\end{align}

The variable distributions are:
\begin{align}
F\equiv{RMS}(v_1^{(100)}) \notag \\
\mathbf{X}_d \sim \mathcal{N}_t(\boldsymbol{\mu}_d,\boldsymbol{\sigma}_d^2), \; \boldsymbol{\sigma}_{d} = 0.02\boldsymbol{\mu}_{d}  \notag\\
\bar{A} \sim U(8,	12) \notag \\
\bar{\omega}_1 \sim U(0.48,0.52) \notag \\
\bar{\omega}_2 \sim U(1.23, 1.27) \notag
\end{align}

Results are listed in Table \ref{tab:opt2_result}. The RMS displacement value is reduced by 39.7\% compared to the reference configuration, whose RMS response is 8.79.

\begin{table}[!h]
	\small
	\centering
	\caption{\small{Problem 2. Optimization results.}}
	\label{tab:opt2_result}
	\begin{tabular}{ccccc}
		\Xhline{2pt}
		\centering
		$\mathbf{\mu}_{\alpha}$& $\mathbf{\mu}_{\epsilon}$& $\mathbf{\mu}_{\eta}$&$\mathbb{E} ({F})$ & Rel. reduction. reference. \\
		\hline
		\centering
		$0.471$ &$2.08$ &$4.72 \times 10^{-4}$ & 5.30 & 39.7\% \\
		\Xhline{2pt}
	\end{tabular}
\end{table}

\begin{figure}[!h]
	\centering
	\includegraphics[width=0.45\linewidth]{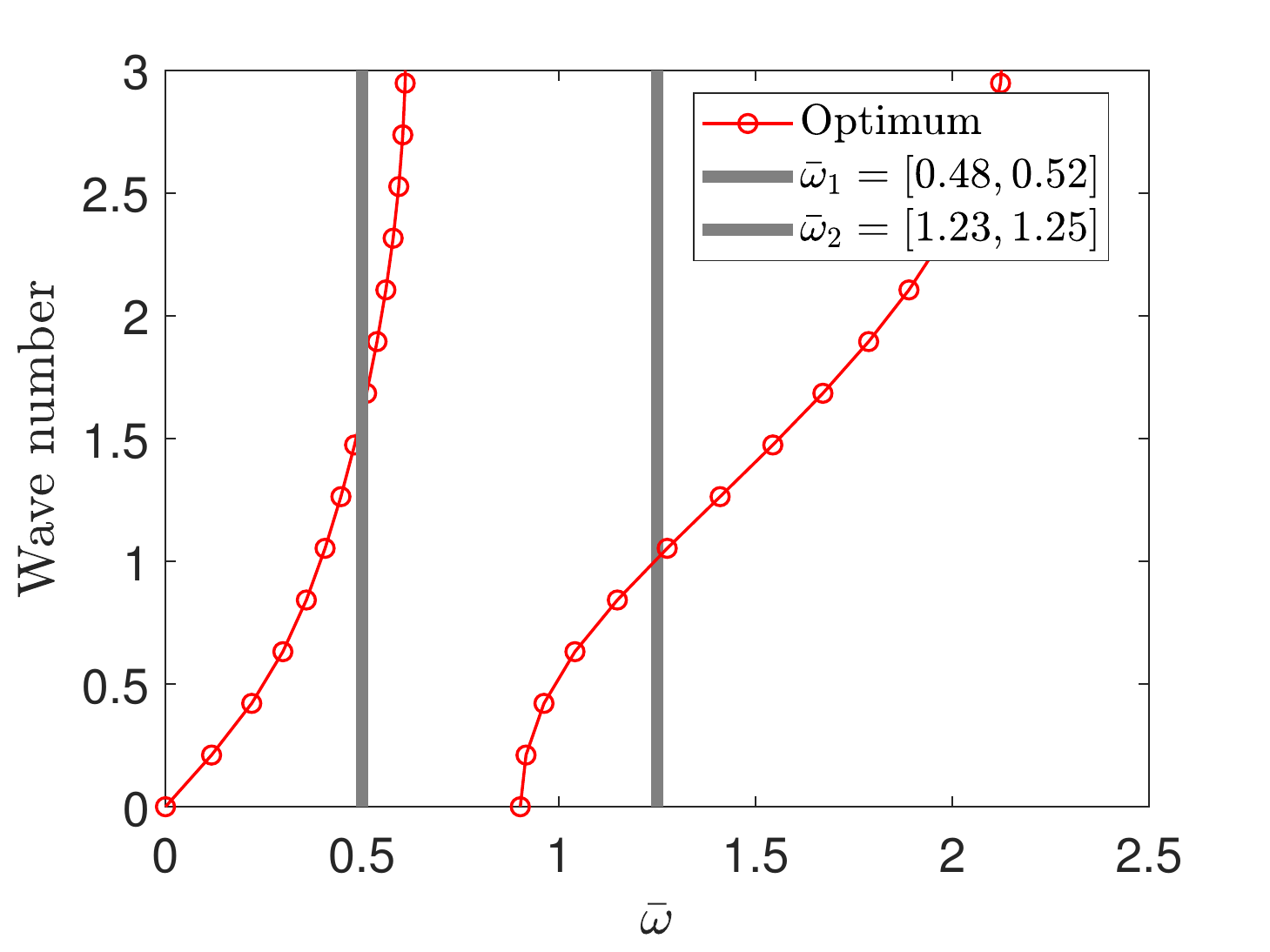}
	\caption{\small{Problem 2. Dispersion diagram of the optimal design {($\alpha = 0.471$, $\epsilon = 2.08$ and $\eta = 4.72\times 10^{-4}$)} for amplitude of $\bar{A}=10$.}}
	\label{fig:optim2_HB}
\end{figure}

The dispersion behavior of the optimally designed chain is illustrated in Figure \ref{fig:optim2_HB}. It shows that none of the excitation frequency ranges falls within the band gap range. It implies that, despite some reduction in amplitude, it is not possible to design a chain, with constant properties over the chain and within the defined boundaries of the design variables, which is able to mitigate both excitation frequency ranges.

\begin{figure}
	\centering
	\begin{subfigure}[b]{0.4\textwidth}
		\includegraphics[width=\linewidth,keepaspectratio]{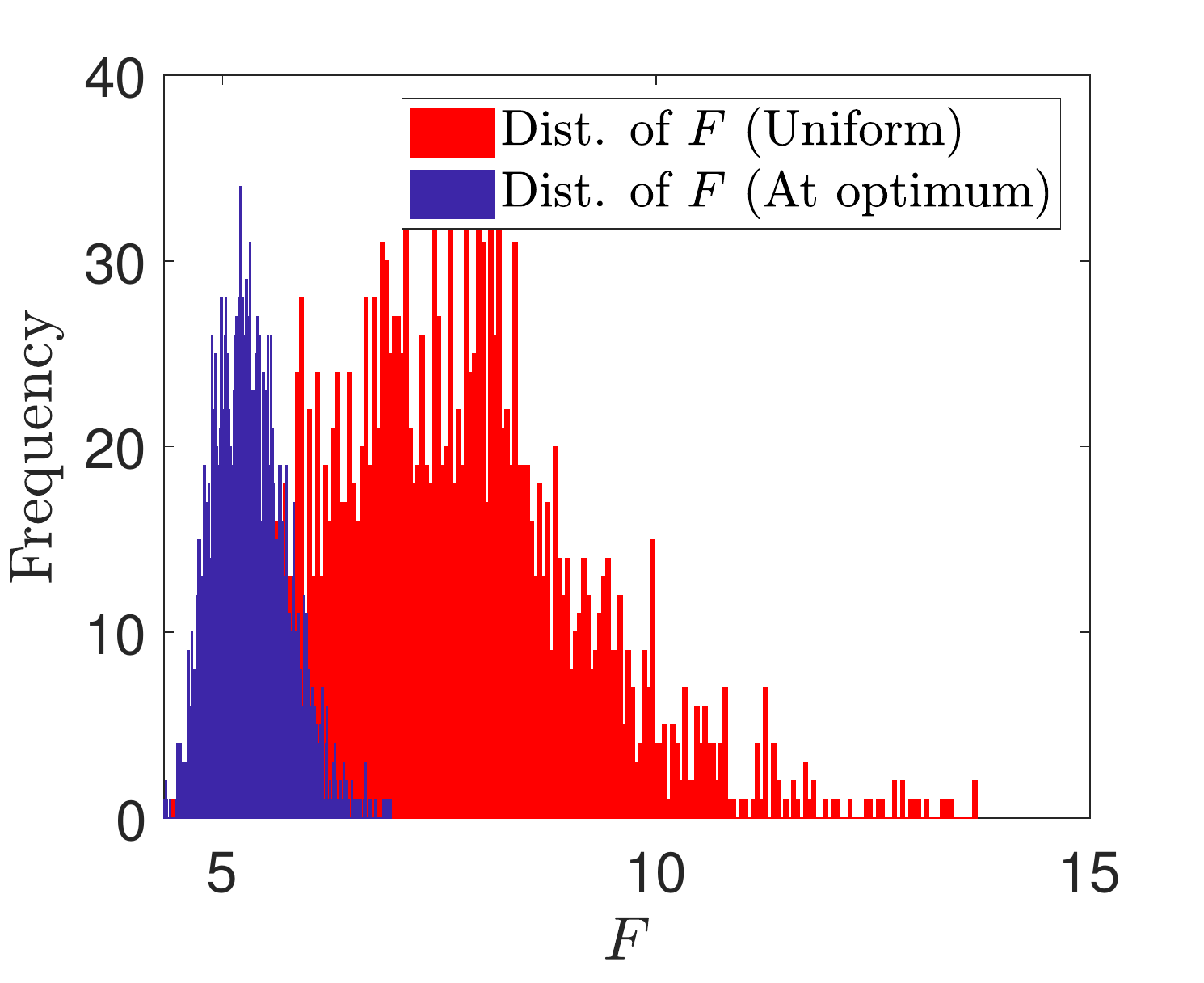}
		\caption{}
		\label{fig:hist_unif_problem2_a}
	\end{subfigure}
	\begin{subfigure}[b]{0.4\textwidth}
		\includegraphics[width=\linewidth,keepaspectratio]{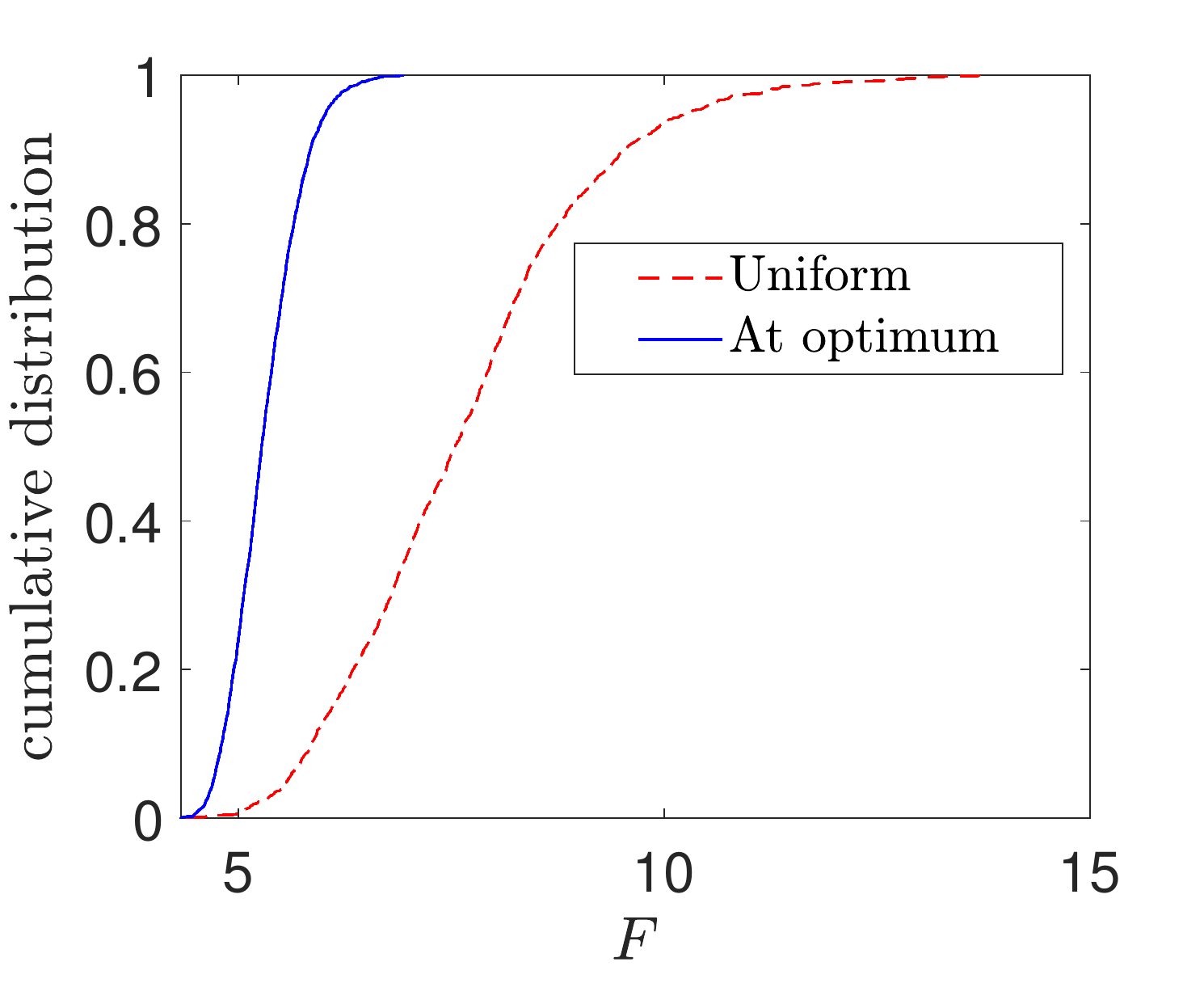}
		\caption{}
		\label{fig:hist_unif_problem2_b}
	\end{subfigure}
	\caption{\small{Problem 2. Comparing the distributions of RMS responses ($F$) around optimum vs. variables uniformly distributed. (a): Histograms of distributions. (b): Cumulative distributions.}}
	\label{fig:hist_unif_problem2}
\end{figure}

The distribution of RMS values around the optimum compared to the distribution of responses assuming uniform distributions for all variables is depicted in Figure \ref{fig:hist_unif_problem2}.  The figure shows that the non-optimal cases with uniform distributions have a much larger spread, thus showing the effect of the optimization which led to a more ``robust" optimum. It should be noted that the difference between the optimal and non-optimal designs is more obvious compared to the previous optimization problem because more design variables are involved.

\subsection{Optimization problem 3. 9 design variables. 3 aleatory variables}
In the third optimization problem, the field function as described in Section \ref{sec:field} is used to represent the properties of the chain. This optimization problem aims at demonstrating the advantage of using spatially varying properties instead of constant properties over the chain (optimization problem 2). The previous optimization problem is repeated except that the design variables follow a field. The optimization variables are the coefficients of the fields, $\mathbf{X}_d=\{\mathbf{a}_{\alpha},\mathbf{a}_{\epsilon}, \mathbf{a}_{\eta}\}$. $\mathbf{a}_{\alpha}$, $\mathbf{a}_{\epsilon}$  and $\mathbf{a}_{\eta}$ are the  field coefficients for $\alpha$, $\epsilon$ and $\eta$, respectively (Eq. \ref{eqn:field}). The optimization problem reads:

\begin{align}
\label{eqn:opt3}
\min_{\boldsymbol{\mu}_{\mathbf{a}}^{\alpha},\boldsymbol{\mu}_{\mathbf{a}}^{\epsilon}, \boldsymbol{\mu}_{\mathbf{a}}^{\eta}} 
& \quad  \mathbb{E} ({F}(\mathbf{a}_{\alpha},\mathbf{a}_{\epsilon}, \mathbf{a}_{\eta}, A,\bar{\omega}_1,\bar{\omega}_2))\\ \notag
s.t.
&\quad \mathbf{a}^{min} + 3 \sigma_{\mathbf{a}} \leq \boldsymbol{\mu}_{\mathbf{a}}\leq \mathbf{a}^{max} - 3 \sigma_{\mathbf{a}} \notag
\end{align}

The ranges of the field coefficients (Table \ref{tab:lin_deter_3omega_field_config}) are chosen so that they span the same ranges as the design variables in the second optimization problem.
\begin{table}[!h]
	\small
	\centering
	\caption{\small{Optimization problem 3. Design variables $\alpha$, $\epsilon$ and $\eta$ are optimized using the field representation. Ranges of field coefficients.}}
	\label{tab:lin_deter_3omega_field_config}
	\begin{tabular}{p{35mm}p{35mm}p{44mm}}
		\Xhline{2pt}
		$\alpha$ & $\epsilon$ & $\eta$   \\
		\hline
		${a_1}\in[10^{-6} ,0.25]$\newline ${a_2}\in[10^{-6},0.25]$\newline ${a_3}\in[10^{-6}, 0.5]$  & ${a_1}\in[0.25,1.25]$\newline ${a_2}\in[10^{-2},1.25]$\newline ${a_3}\in[10^{-2},0.5]$  &${a_1}\in[10^{-6} ,2.5 \times10^{-4}]$\newline ${a_2}\in[10^{-6} ,2.5 \times10^{-4}]$\newline ${a_3}\in[10^{-6},0.5]$  \\
		\Xhline{2pt}
	\end{tabular}
\end{table}

The optimization results which are listed in Table \ref{tab:opt3_result} shows a reduction of 85.6\% compared to the optimization problem 2 which had constant properties over the chain (5.30 to 0.77 RMS displacement). Therefore the field enabled a much drastic reduction of RMS displacement compared to the constant property case. The reduction is 91.2\% compared to the reference design.  
\begin{table}[!h]
	\small
	\centering
	\caption{\small{Problem 3. Design variables $\alpha$, $\epsilon$ and $\eta$ are optimized using the field representation. Optimization results.}}
	\label{tab:opt3_result}
	\begin{tabular}{p{18mm}p{18mm}p{27mm}ccc}
		\Xhline{2pt}
		\centering
		$\boldsymbol{\mu}_{\mathbf{a}}^{\alpha}$& $\boldsymbol{\mu}_{\mathbf{a}}^{\epsilon}$& $\boldsymbol{\mu}_{\mathbf{a}}^{\eta}$&$\mathbb{E} ({F})$ &Rel red pb2 &Rel red ref  \\
		\hline
		$a_1=0.235$ \newline $a_2=0.235$ \newline $a_3=0.336$ &$a_1=0.266$ \newline $a_2=1.18$ \newline $a_3=0.377$&$a_1=1.10\times10^{-6}$ \newline $a_2=1.10\times10^{-6}$ \newline $a_3=1.10\times10^{-6}$ & 0.770 & 85.5 \% &91.2\% \\
		\Xhline{2pt}
	\end{tabular}
\end{table}

\begin{figure}[!h]
	\centering
	\includegraphics[width=0.6\linewidth]{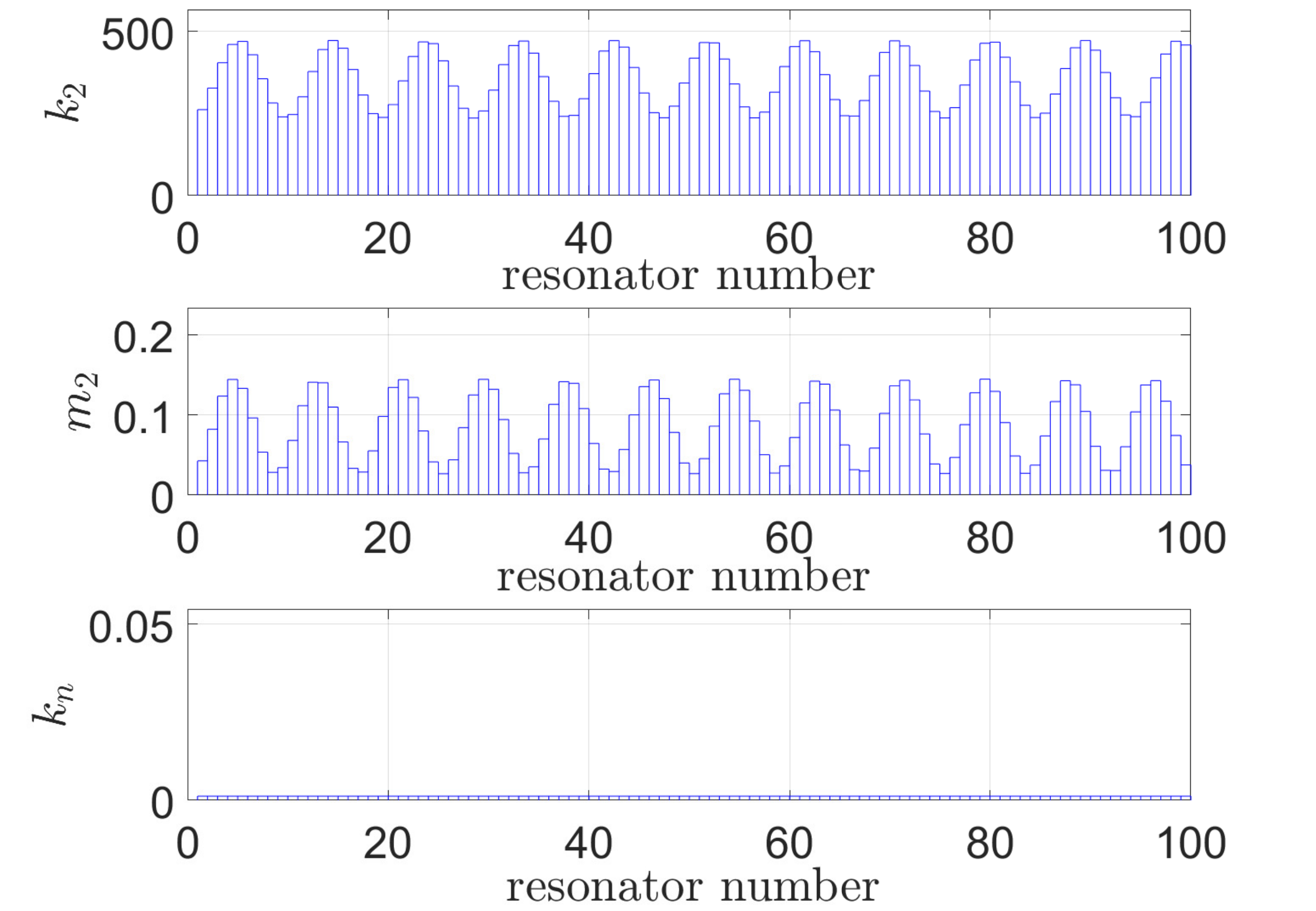}
	\caption{\small{Problem 3. Design variables $\alpha$, $\epsilon$ and $\eta$ are defined using the field representation. Variation of properties over 100 units of the chain (in dimensional form with $k_1 = 1000$ and $m_1 = 0.1$) based on the optimal field coefficients.}}
	\label{fig:optim3_config}
\end{figure}

The optimal distribution of the properties of the internal resonators (in a dimensional form) is depicted in Figure \ref{fig:optim3_config}. It shows that due to the implementation of the field, the properties vary periodically over the chain, but with different periods. It is also noteworthy that in this case, the nonlinearity is not needed and its coefficients are reduced to the lower bound.

\begin{figure}[!h]
	\centering
	\includegraphics[width=0.4\linewidth]{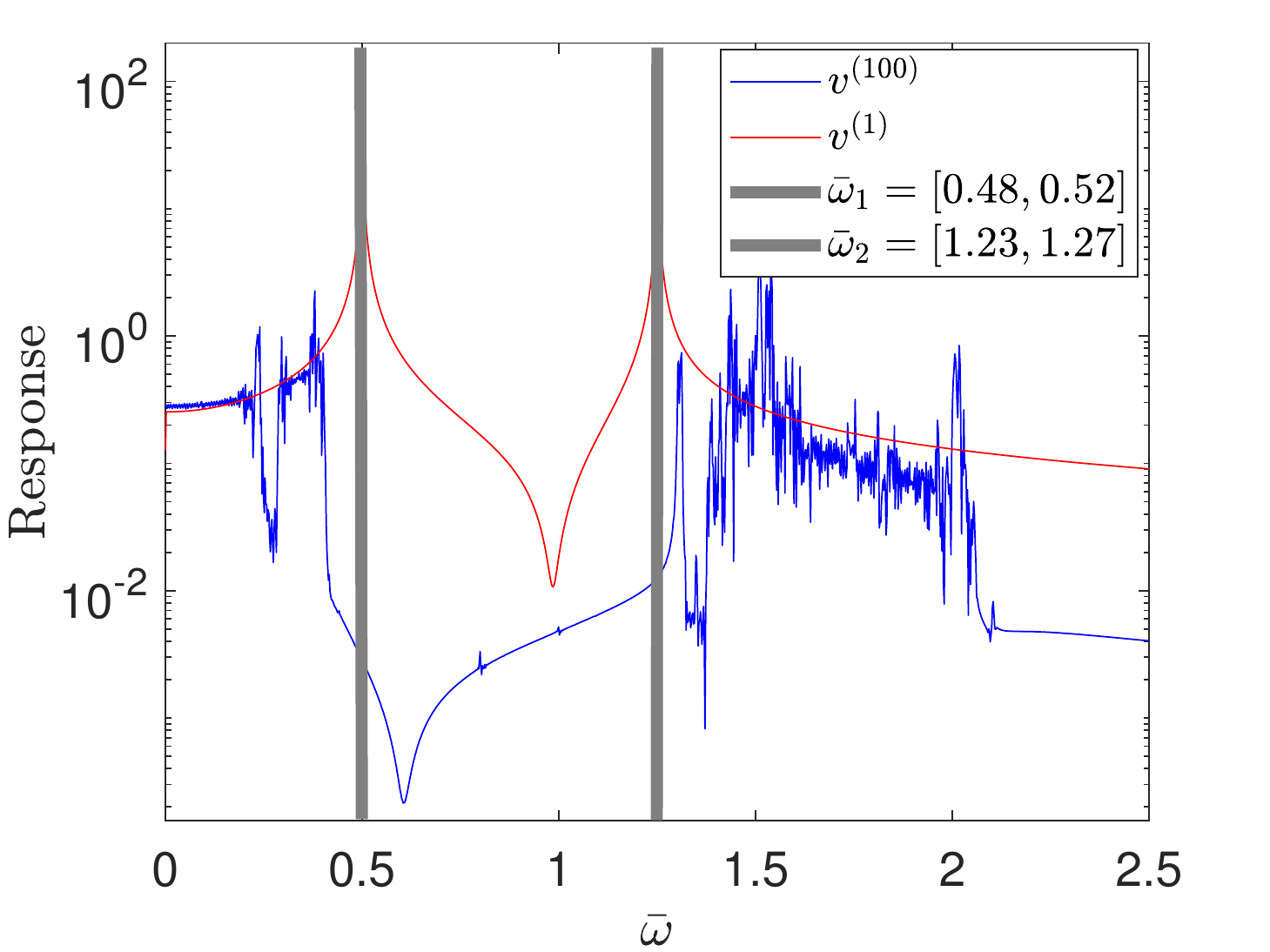}
	\caption{\small{Problem 3. {Frequency responses of the main resonators in unit 1, $v^{(1)}_1$, and unit 100 $v^{(100)}_1$. The use of the field has enabled an expansion of the band gap to include both excitations (with frequencies of $\bar{\omega}_1=0.50$ and $\bar{\omega}_2=1.25$ and amplitude of $\bar{A}=10$) that are applied at unit 1.}}}
	\label{fig:optim3_fft_2_log}
\end{figure}

In order to analyze the results and understand why the field-based design led to a much larger reduction, the frequency response of unit 100, {$v^{(100)}_1$}, of the optimal chain being subjected to a harmonic displacement (with two frequencies of $\bar{\omega}_1=0.50$ and $\bar{\omega}_2=1.25$ and amplitude of $\bar{A}=10$) at unit 1 is depicted in Figure \ref{fig:optim3_fft_2_log}. These values correspond to the means of the aleatory variables' distributions. {For comparison, the frequency response of the main resonator in unit 1, $v^{(1)}_1$, where the excitations are applied is also provided.} The figure clearly shows that the band gap frequency range has been extended to accommodate both frequencies. This was made possible by the introduction of the field description which relaxed the optimization problem compared to the case where the properties are constant over the chain (i.e., problem 2).

\begin{figure}[!h]
	\centering
	\begin{subfigure}[b]{0.4\textwidth}
		\includegraphics[width=\linewidth,keepaspectratio]{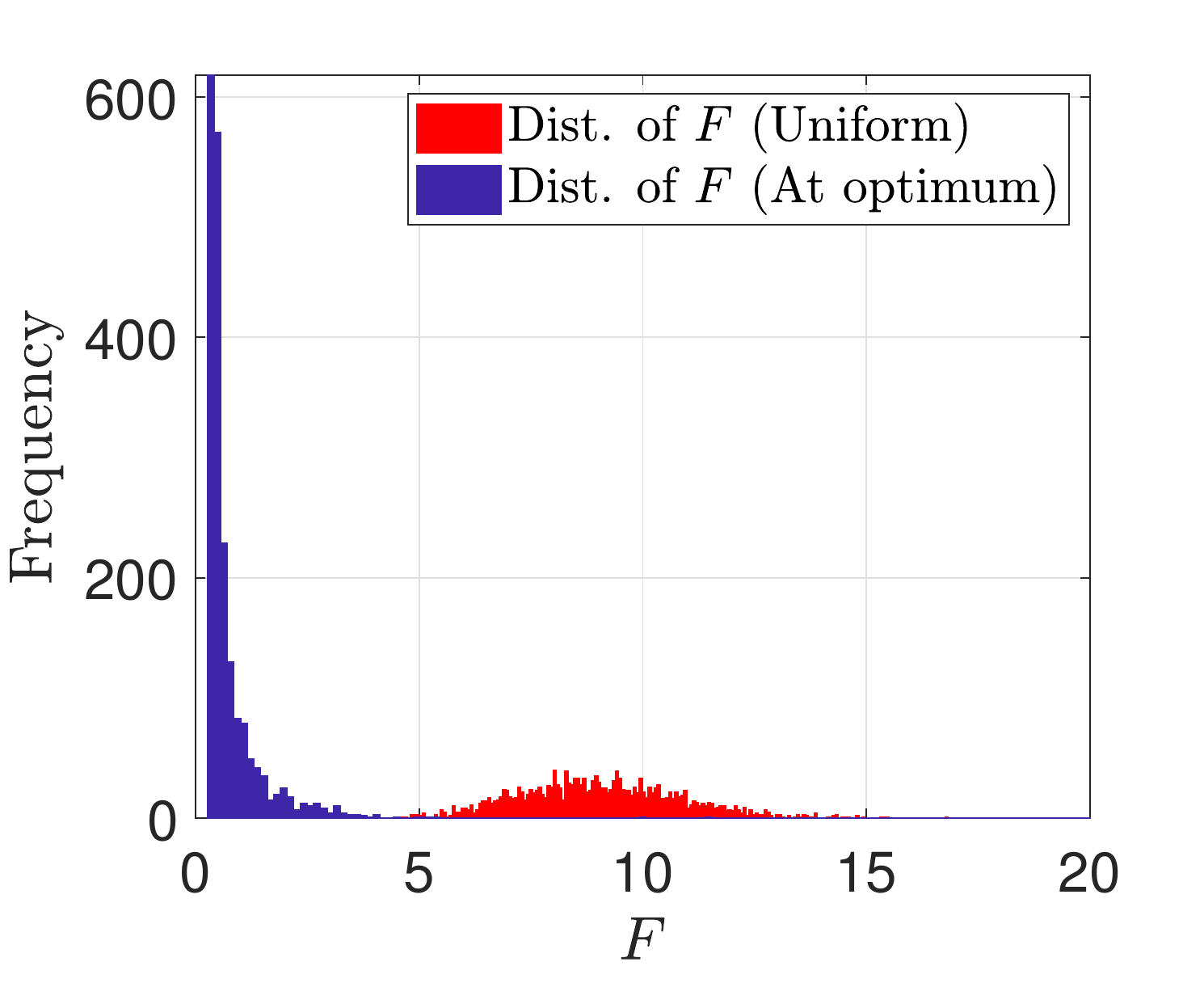}
		\caption{}
		\label{fig:hist_unif_problem3_a}
	\end{subfigure}
	\begin{subfigure}[b]{0.4\textwidth}
		\includegraphics[width=\linewidth,keepaspectratio]{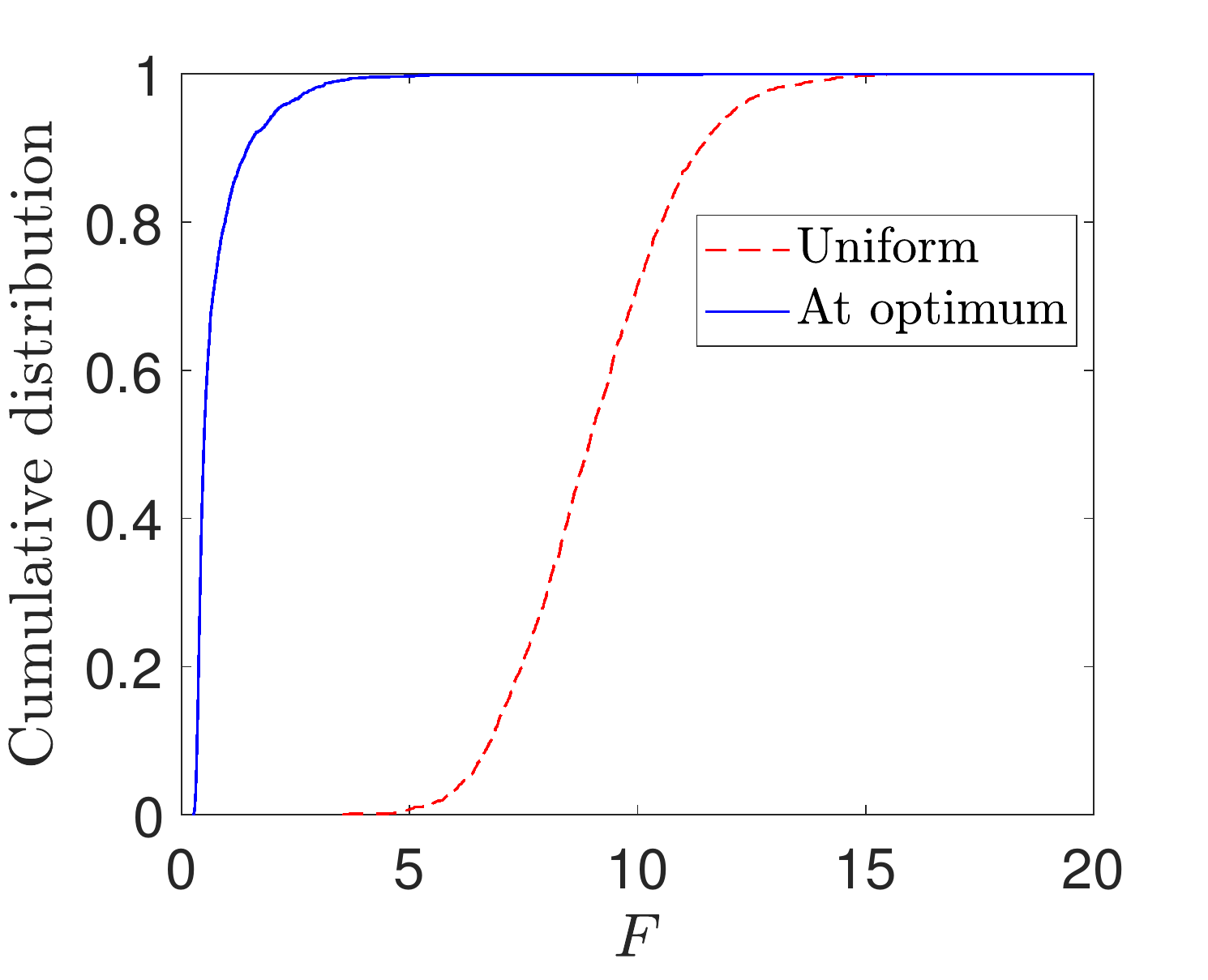}
		\caption{}
		\label{fig:hist_unif_problem3_b}
	\end{subfigure}
	\caption{\small{Problem 3. Comparing the distributions of RMS responses ($F$) around optimum vs. variables uniformly distributed. (a): Histograms of distributions. (b): Cumulative distributions.}}
	\label{fig:hist_unif_problem3}
\end{figure}

Regarding the robustness of the optimization results, Figure \ref{fig:hist_unif_problem3} shows the substantially better performance of the optimal designs compared to the non-optimal cases.

\section{Conclusion}
In this work, the stochastic optimization of a one-dimensional chain of nonlinear resonators was presented for the purpose of vibration mitigation. The optimization approach is capable of tackling response discontinuities due to the presence of nonlinearities. It also accounts for uncertainties in design variables and loading conditions through the use of machine learning techniques such as clustering, SVM and Kriging. In addition, to make the optimization problem scalable, this method makes use of a field representation to define properties such as mass and stiffness over the chain. The use of a field representation enables one to optimize the properties of all resonators using a few variables thus substantially reducing the number of optimization variables. As a demonstrative example of the methodology, it was shown, through a twelve-dimensional problem, that the field description enables an extension of the band gap to minimize the chain's response. Such an extension was limited when constant properties over the chain were used.

Based on these initial results, the next steps of this research will focus on optimally designing 2D and 3D structures, including considering the influence of damping with the possibility of targeted energy transfers between the main and internal resonators. In addition, the use of multiple basis functions in the field representation will be investigated to further improve the achieved optimal design.

\section{Appendix}
\setcounter{table}{0}
\renewcommand{\thetable}{A\arabic{table}}
\setcounter{figure}{0}
\renewcommand{\thefigure}{A\arabic{figure}}
The accuracy of the optimization results, obtained using surrogate models must be validated based on the actual models' responses. For this purpose, the ``true" $\mathbb{E}(F)$ is also computed using a Monte-Carlo simulation using the actual model.  We must also account for the variability due to Monte-Carlo sampling in estimating $\mathbb{E}(F)$. To this end, a histogram of approximations of $\mathbb{E}(F)$ is constructed based on several Monte-Carlo runs. This histogram also enables a correction of the predicted mean value by taking the mean value of the histogram.

\subsection{Problem 1}
\begin{figure}[!h]
	\centering
	\includegraphics[width=0.45\linewidth]{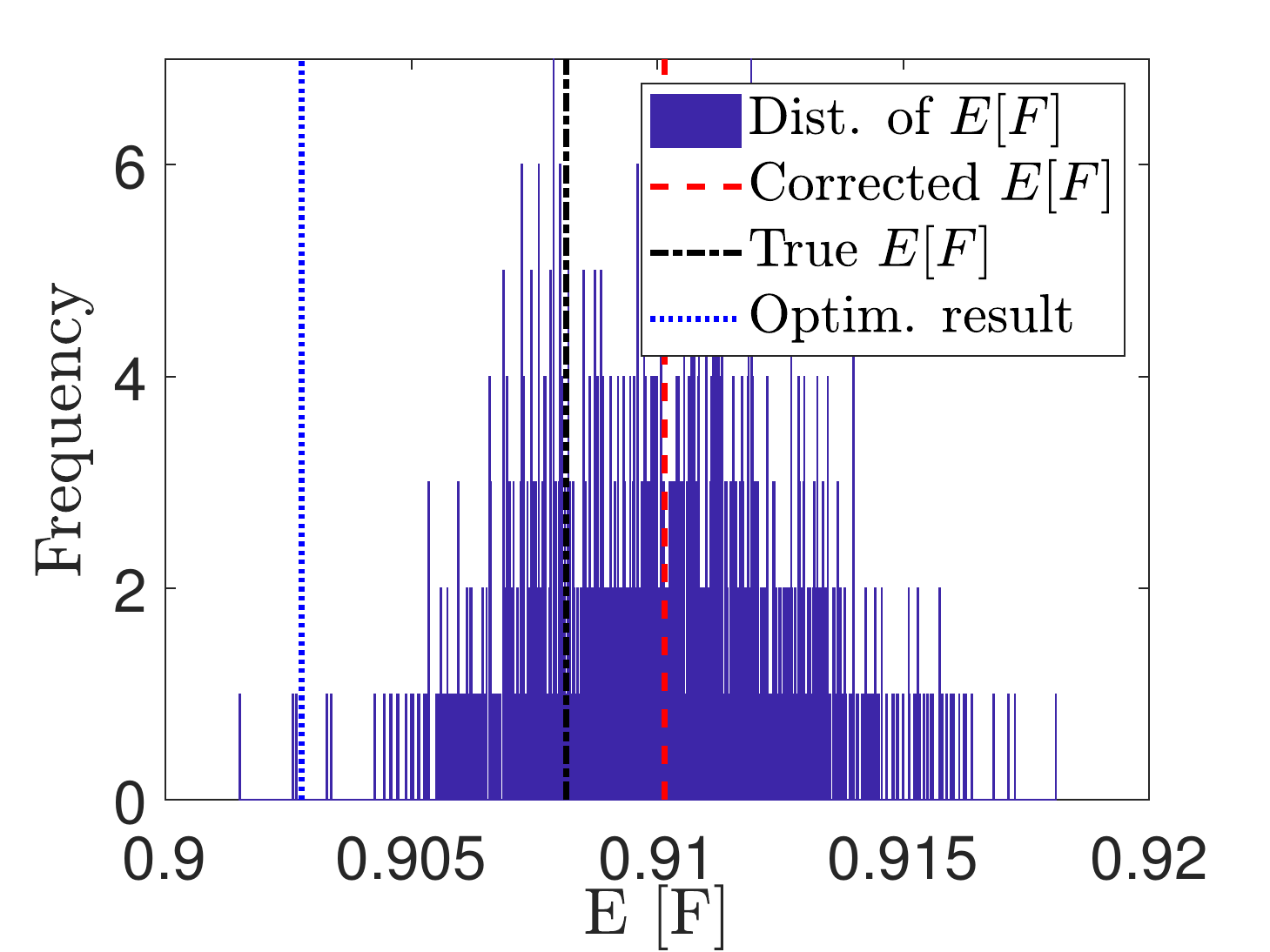}
	\caption{\small{Problem 1. Accuracy of the optimization results. Variability of the objective function at the optimum  due to Monte-Carlo sampling.}}
	\label{fig:optim1_hist}
\end{figure}

The relative error between surrogate-based $\mathbb{E}(F)$ and true $\mathbb{E}(F)$ at the optimum is provided in Table \ref{tab:opt1_result_ap}. The histogram in Figure \ref{fig:optim1_hist} shows that the true $\mathbb{E} ({F})$ falls within the distribution of the surrogate-based $\mathbb{E} ({F})$ at the optimum with relative error of 0.2\%. Note that since PSO is a population-based optimization method and returns the best in a population, the optimization result generally falls at the lower tail of the distribution plotted in Figure \ref{fig:optim1_hist}. The optimum value is corrected as being the mean of this distribution and is reported in the table. Note that the results reported in Section \ref{sec:results} are based on the corrected value of $\mathbb{E} ({F})$. \\

\begin{table}[!h]
	\small
	\centering
	\caption{\small{Problem 1. Checking the accuracy of the optimization results.}}
	\label{tab:opt1_result_ap}
	\begin{tabular}{cccc}
		\Xhline{2pt}
		\centering
		$\mathbb{E} ({F})$ &Corrected $\mathbb{E} ({F})$ & True $\mathbb{E} ({F})$& Rel. error \\
		\hline
		\centering
		0.903&0.910 &0.908 & 0.2\%\\
		\Xhline{2pt}
	\end{tabular}
\end{table}

\begin{figure}[!h]
	\centering
	\includegraphics[width=0.35\linewidth]{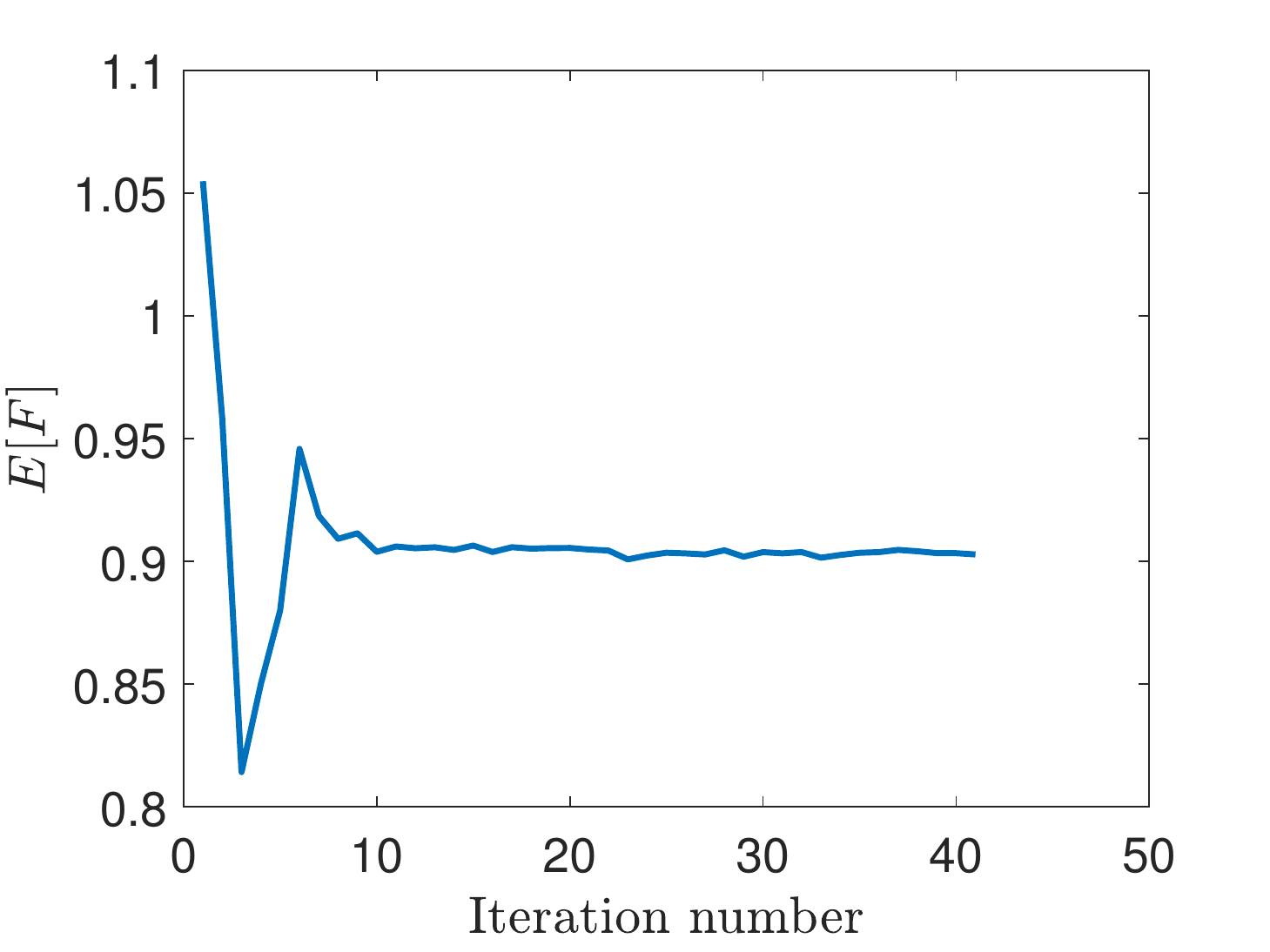}
	\caption{\small{Problem 1. Evolution of the objective function.}}
	\label{fig:optim1_iter}
\end{figure}

The evolution of the objective function is depicted in Figure \ref{fig:optim1_iter}. It can be seen that the optimization converges in about 10 iterations which corresponds to 20 evaluations {of the actual computational model} in addition to the 10 samples used in the initial design of experiments. {It is noteworthy that 2 new samples are generated in each iteration to refine the Kriging and SVM metamodels.}

\subsection{Problem 2}

\begin{figure}[!h]
	\centering
	\includegraphics[width=0.45\linewidth]{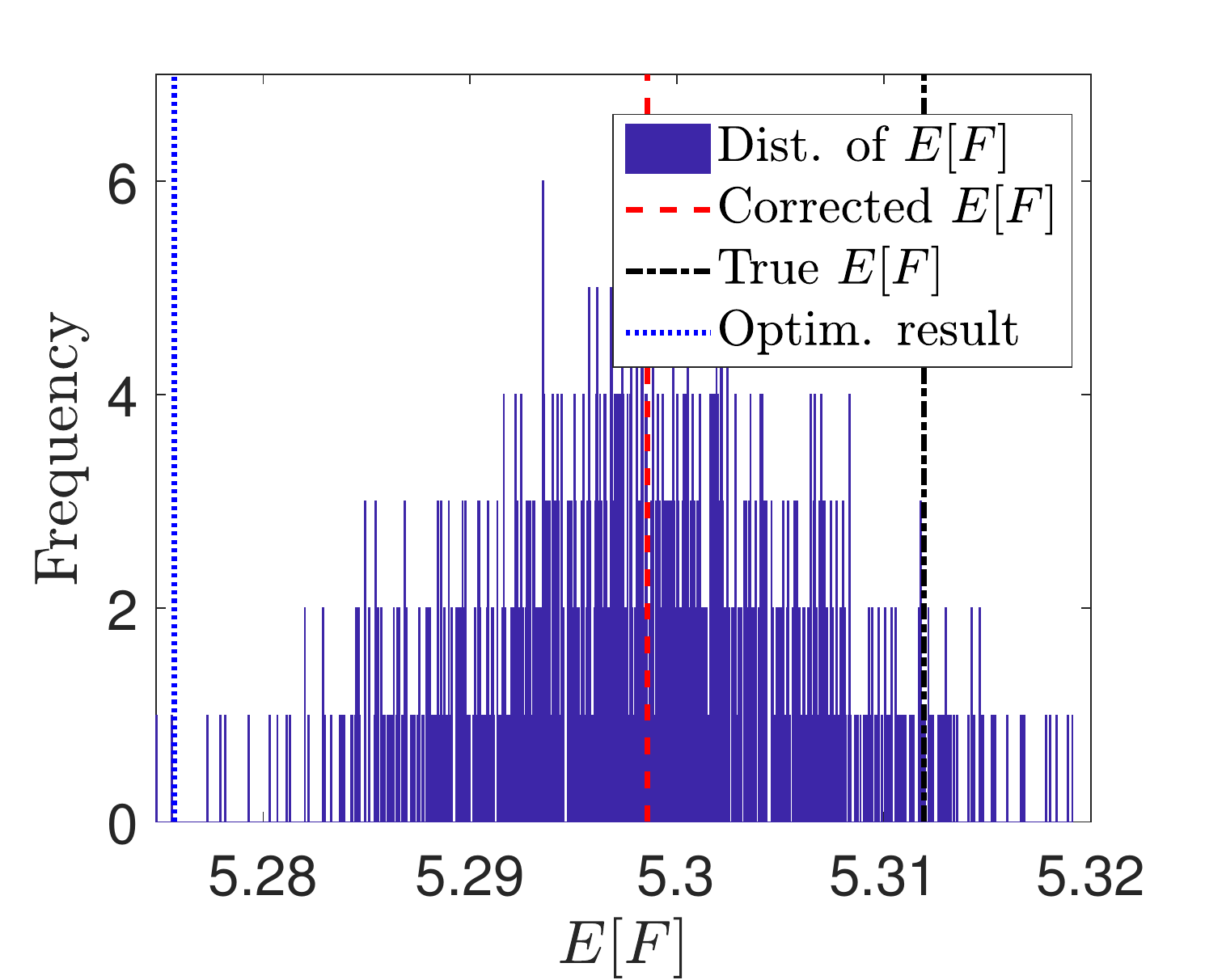}
	\caption{\small{Problem 2. Accuracy of the optimization results. Variability of the objective function at the optimum  due to Monte-Carlo sampling.}}
	\label{fig:optim2_hist}
\end{figure}

\begin{figure}[!h]
	\centering
	\includegraphics[width=0.35\linewidth]{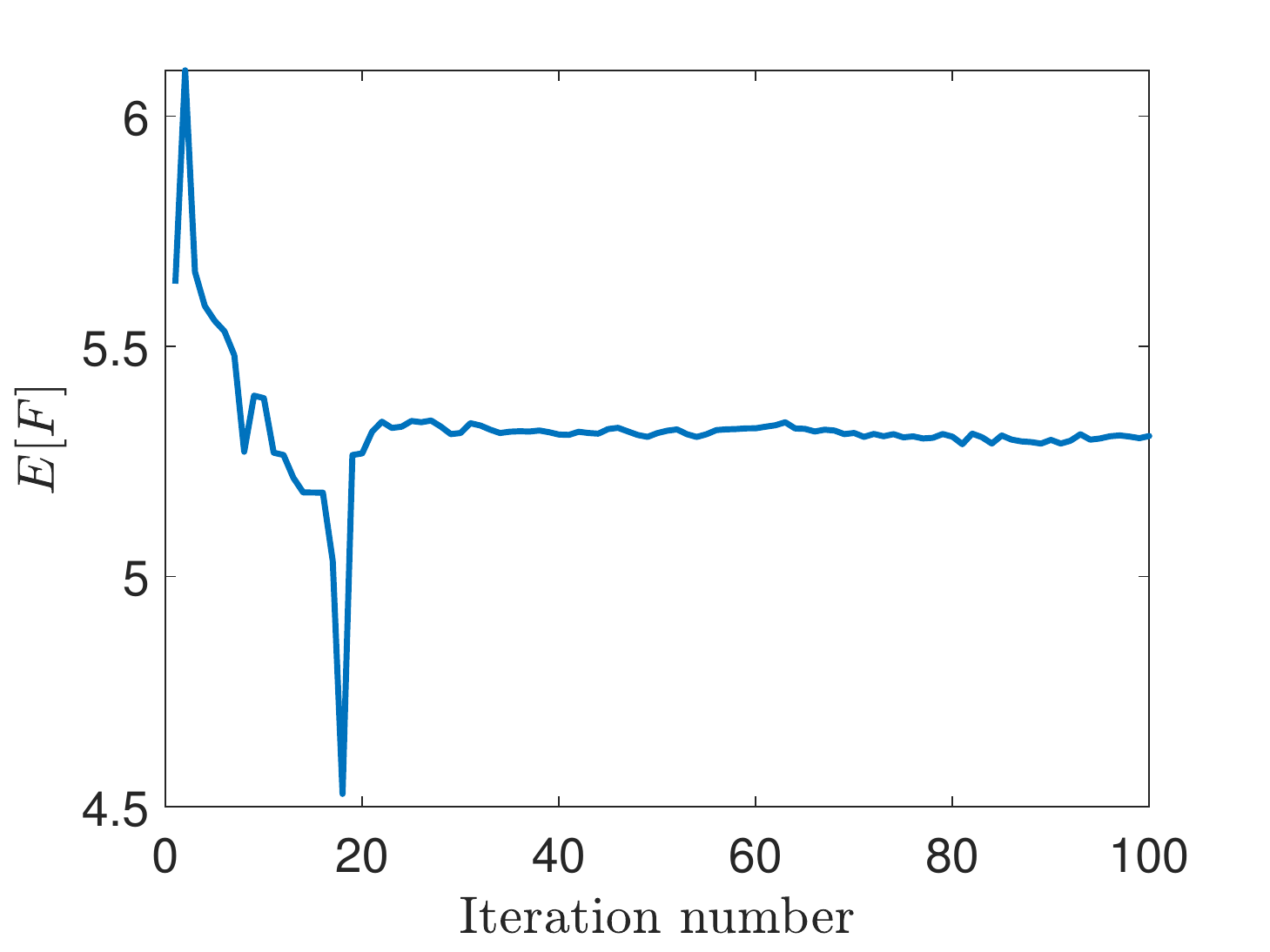}
	\caption{\small{Problem 2. Evolution of the objective function.}}
	\label{fig:optim2_iter}
\end{figure}

A similar study was carried out for the second problem. The relative error is provided in Table \ref{tab:opt2_result_ap}. Figure \ref{fig:optim2_hist} shows that the true $\mathbb{E}(F)$ falls within the distribution of $\mathbb{E}(F)$ at the optimum with relative error of 0.19\%. The evolution of the objective function is depicted in Figure \ref{fig:optim2_iter}. 
An initial design of experiments of size 75 is used in this optimization.
\begin{table}[!h]
	\small
	\centering
	\caption{\small{Problem 2. Checking the accuracy of the optimization results.}}
	\label{tab:opt2_result_ap}
	\begin{tabular}{cccc}
		\Xhline{2pt}
		\centering
		$\mathbb{E} ({F})$ &Corrected $\mathbb{E} ({F})$ & True $\mathbb{E} ({F})$& Rel. error \\
		\hline
		\centering
		5.28&5.30 & 5.31 & 0.19\%\\
		\Xhline{2pt}
	\end{tabular}
\end{table}

\subsection{Problem 3}

\begin{figure}[!h]
	\centering
	\includegraphics[width=0.45\linewidth]{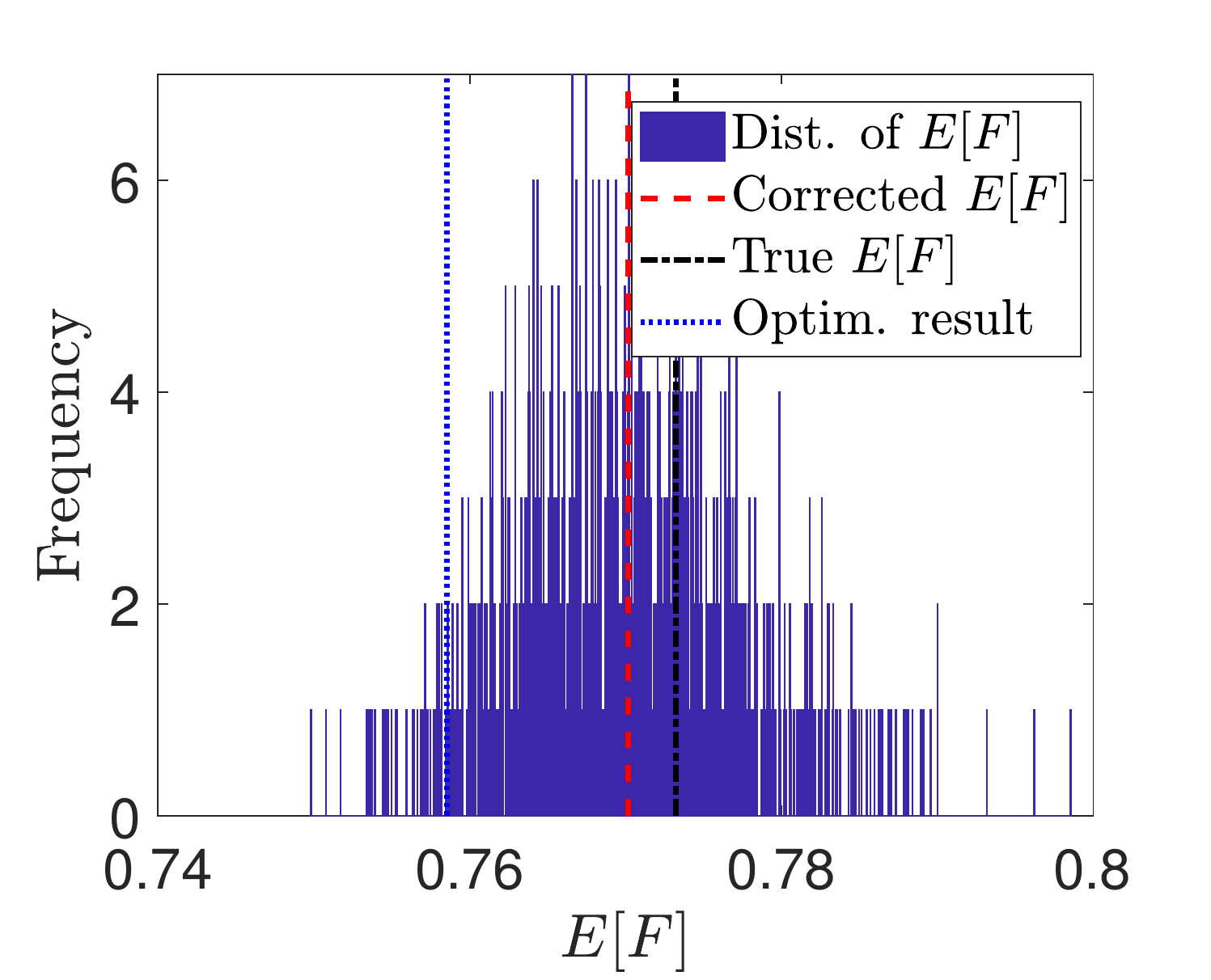}
	\caption{\small{Problem 3. Accuracy of the optimization results. Variability of the objective function at the optimum  due to Monte-Carlo sampling.}}
	\label{fig:optim3_hist}
\end{figure}

\begin{figure}[!h]
	\centering
	\includegraphics[width=0.35\linewidth]{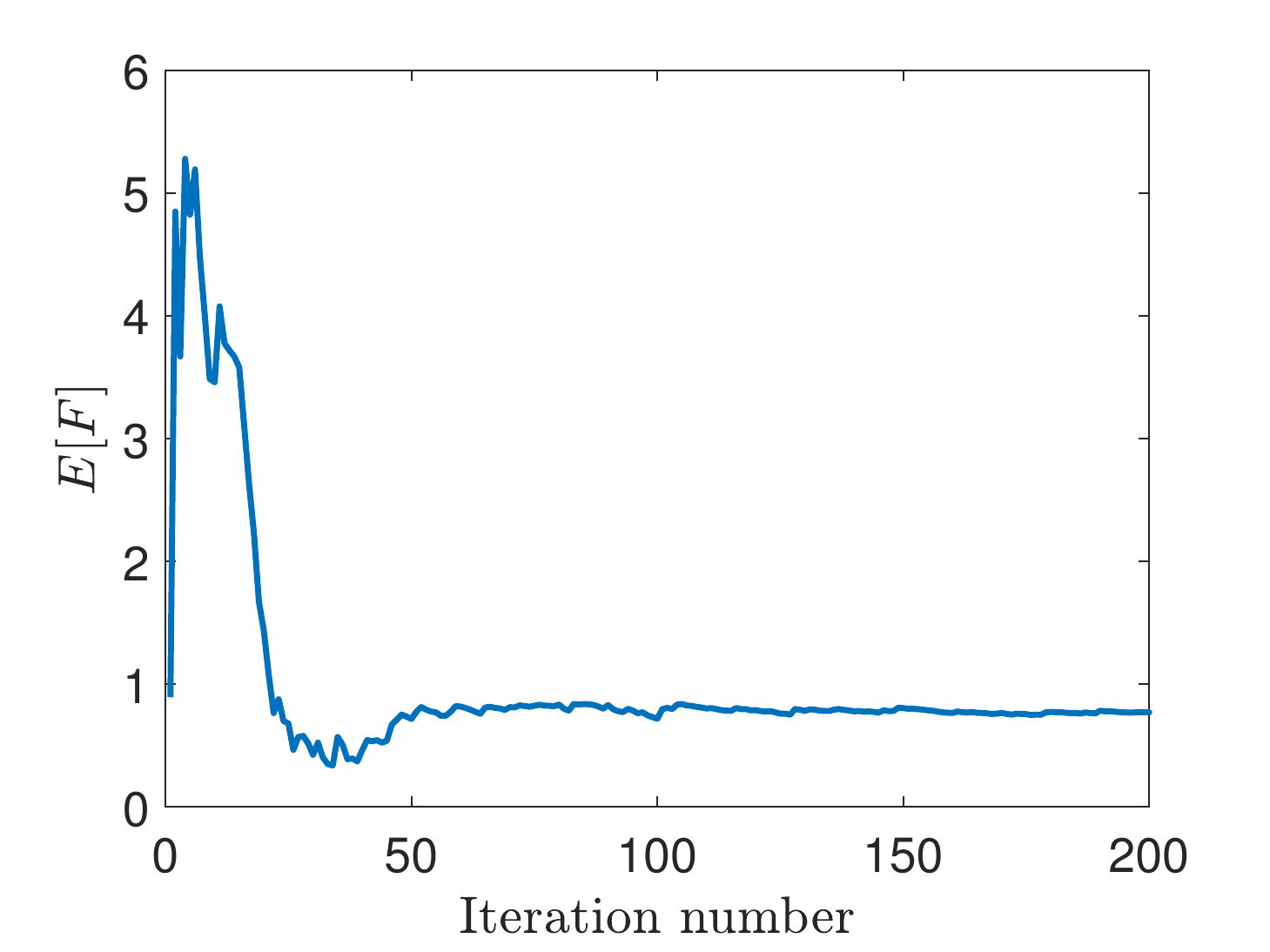}
	\caption{\small{Problem 3. Evolution of the objective function.}}
	\label{fig:optim3_iter}
\end{figure}

Similar to optimization problems 1 and 2, the accuracy of the results is checked, and the results are provided in Table \ref{tab:opt3_result_ap} and Figure \ref{fig:optim3_hist}. The evolution of the objective function is provided in Figure \ref{fig:optim3_iter}. 
175 samples are used in the initial design of experiments.
\begin{table}[!h]
	\small
	\centering
	\caption{\small{Problem 3. Checking the accuracy of the optimization results.}}
	\label{tab:opt3_result_ap}
	\begin{tabular}{cccc}
		\Xhline{2pt}
		\centering
		$\mathbb{E} ({F})$&Corrected $\mathbb{E} ({F})$ & True $\mathbb{E} ({F})$& Rel. error \\
		\hline
		\centering
		0.759 & 0.770 &0.773 & 0.39\%\\
		\Xhline{2pt}
	\end{tabular}
\end{table}

\section{References}
\vspace{-0.2cm}
\bibliographystyle{unsrt}
\bibliography{mybiblio}

\biboptions{sort&compress}
\end{document}